\def\**{$*\!*\!*$}


\font\bit=cmssi12 at 12truept
\font\tenmsy=msym10
\textfont8=\tenmsy
\mathchardef\ssm="7872

\input psfig
\mathsurround = 2pt
\abovedisplayskip=6pt
\belowdisplayskip=6pt

\def \QP{\narrower\medskip\noindent}
\def \QED {\rlap{$\sqcup$}$\sqcap$\smallskip}
\def\Z{{\bf Z}}
\def\R{{\bf R}}

\def\ref{\hangindent=1pc \hangafter=1 \noindent}
\def\mod{{\rm mod\;}}
\def\sgn{{\rm sgn}}


\centerline{\bf The Fibonacci Unimodal Map.}\medskip

\centerline{Mikhail Lyubich and John Milnor}\smallskip
\centerline{Mathematics Department and IMS, SUNY Stony Brook}\bigskip
\line{\hfill {\it To Yuri Lyubich on his 60th birthday.}\hskip 1in}\bigskip

\centerline{\bf \S1. Introduction.}\smallskip

The Fibonacci recurrence of the critical orbit appeared in the work of Branner
and Hubbard on complex cubic polynomials [BH, \S12], and in Yoccoz's work [Y]
on quadratic ones,  as the ``worst'' pattern of recurrence. On the other hand,
a real quadratic Fibonacci map $f$ was suggested by Hofbauer and Keller [HK]
as a possible candidate for a map having a ``wild" attractor (because the $\omega$-limit
set of the critical point possesses all known topological
properties of wild attractors [BL2] ).  Also, Shibayama [Sh] has described
this real Fibonacci map as the limit of a sequence of quadratic maps with
attracting orbit whose period is a Fibonacci number.

This paper will study topological,
geometrical and measure-theoretical properties of the real Fibonacci map.
Our goal was to figure out if this type of recurrence really gives any
pathological examples and to compare it with the infinitely renormalizable
patterns of recurrence studied by Sullivan [S].
It turns out that the situation can be understood completely and is of quite
regular nature. 
In particular, any Fibonacci map (with negative Schwarzian and non-degenerate
critical point) has an absolutely continuous 
invariant measure (so, we deal with a ``regular'' type of chaotic dynamics).
It turns out also that geometrical properties of the closure of the critical
orbit are quite different from those of the Feigenbaum map: its Hausdorff 
dimension is equal to zero and its geometry is not rigid but depends on one 
parameter. 

 Branner and Hubbard
introduce the concept of a {\bit tableau\/} in order to
describe recurrence of critical orbits. Their ``Fibonacci tableau''
is a basic example, which corresponds to
one particularly close and regular pattern of recurrence.
If a complex quadratic map $z\mapsto z^2+c$ realizes this
Fibonacci tableau, then the orbit
$$	0=z_0\mapsto z_1\mapsto z_2\mapsto\cdots $$
of the critical point returns closer to
zero (in a certain invariant sense)
after each Fibonacci number of iterations. In the real case, it follows that
$$	|z_1|\,>\,|z_2|\,>\,|z_3|\,>\,|z_5|\,>\,|z_8|\,>\,|z_{13}|\,>\,
	\cdots\,.$$
In \S2 we will prove that a real quadratic map is uniquely defined by the
last property. More precisely we prove the following.
We denote the Fibonacci numbers by
$$ u(1)=1\,,\;u(2)=2\,,\;\ldots\,,\qquad{\rm with}\quad
u(n+1)=u(n)+u(n-1)\,.$$

{\QP{\bf Theorem 1.1.} \it There is one and only one real quadratic map
of the form\break $\;f_c(x)= x^2+c\;$
with the property that the critical orbit $\;0=x_0\mapsto x_1\mapsto\cdots\;$ 
has  closest recurrence at the
Fibonacci values, so that $\;|x_1|>|x_2|>|x_3|>|x_5|>\cdots\,$, with $x_4<0
\;$.
\footnote\dag{\rm We conjecture that this condition on $x_4$ is automatically
satisfied.}
The kneading invariant for this uniquely defined map $f_c$
can be described by the conditions that
$$\eqalign{
	x_{u(n)}<0\quad &{\rm for}\quad n\equiv 0\,,\,1\quad\mod 4\cr
	x_{u(n)}>0\quad &{\rm for}\quad n\equiv 2\,,\,3\quad\mod 4\;,\cr}$$
and that
$$	\sgn(x_i)\;=\;\sgn(x_{i-u(n)})\qquad{\rm for}\quad u(n)<i<u(n+1)\,.$$
\medskip}

\noindent In fact numerical computation shows that
$\quad c\; =\; -1.8705286321646448888906\cdots \;. $\break
The associated topological entropy is $\;\;h\;=\;\log\,1.7292119317\cdots\;$.
\bigskip

For a fairly general unimodal map $f$
with this same kneading data, we prove the following. Let ${\cal O}=\{x_0\,,\,
x_1\,,\,\ldots\}\subset\R$ be the critical orbit.\medskip

{\QP{\bf Theorem 1.2.} \it If $f$ is $C^2$-
smooth with non-flat critical point, and with kneading data as above,
then:\smallskip


\item{\rm 1.} The closure $\overline{\cal O}$ of the
critical orbit is a Cantor set, 
with the $x_i\;,
\;\;i\ge 1,$ as the end points of the complementary intervals.\smallskip

\item{\rm 2.} The map $f$ from $\overline{\cal O}$ onto itself is one-to-one
except that the critical point has two pre-images. This map
$f|\overline{\cal O}$ is minimal, and is uniquely ergodic with entropy zero.
It is semi-conjugate to the golden rotation
$$ 	t\;\mapsto t-(\sqrt 5-1)/2\qquad (\mod 1)$$
of the circle $\R/\Z$.\medskip}

\noindent
The proof, in \S3, will give an explicit description of the ordering of this
critical orbit closure. It will also show that it is canonically homeomorphic
to the set of all infinite sequences $(a_1\,,\,a_2\,,\,\ldots\,)$ of zeros
and ones with no two consecutive ones, or to the set of all finite or infinite
``Fibonacci sums''. (Compare  2.3 and 3.3.)\smallskip

{\QP{\bf Theorem 1.3.} \it If $f$ is $C^2$-smooth with non-degenerate critical
point then:\smallskip

\item{\rm 1.} The ratio of  $x_{u(n)}$ to $x_{u(n-1)}$  decreases
exponentially, with
$$ \lambda_n\equiv |x_{u(n)}|/|x_{u(n-1)}|\;\sim\; a/2^{n/ 3} 
\qquad{\rm as}\qquad n\to\infty\; $$ for some constant $a>0$.

\item{\rm 2.} The critical orbit closure $\overline{\cal O}$ has
Hausdorff dimension zero and the Liapunov exponent at the critical value
is equal to zero.\smallskip

\item{\rm 3.} Any two Fibonacci maps  
with the same parameter $a$ are smoothly conjugate on 
$\overline{\cal O}$ .\smallskip

\item{\rm 4.} If the Schwarzian derivative is negative, then
$f$ has a unique absolutely continuous invariant measure,
with support equal to the entire closed interval
$\,[x_1\,,\,x_2]\,$, and with positive entropy .\medskip }

{\bf Remark 1.} Uniqueness and other properties of an absolutely continuous 
invariant measure hold automatically (see [BL2].) Existence we will derive
from the Nowicki-van Strien ``series" condition [NvS].

{\bf Remark 2.}  
Unlike the Feigenbaum map, the geometry of $\overline{\cal O}$
goes down to zero under renormalization, and
is not rigid but depends on the parameter $a$.
(We can effectively  vary this parameter).

{\bf Remark 3.} It is essential here that the critical point be
non-degenerate ($f''(0)\ne 0$). We hope to show in a later paper that, for
example, a Fibonacci map of the form $f(x)=x^4+c$ has completely different
behavior, with bounded geometry and with no absolutely continuous
invariant measure.\smallskip

Let us describe the structure of the proof of the last theorem, which
is somewhat complicated.
In \S 4 we get some a priori bounds on the ratios $\lambda_n$. In \S 5 we
prove the Theorem assuming that $\inf\lambda_n=0$. In order to verify this
assumption we introduce in \S 6 an appropriate notion of renormalization so that
infinitely renormalizable maps are exactly Fibonacci ones. Applying Sullivan's
ideas [S] to our case  we prove that if geometry of $\overline{\cal O}$ 
is bounded  from below
then there is a sequence of renormalizations converging to a map which can
be analytically continued in a quite big domain of the complex plane.
   
  In \S 7 we discuss polynomial-like maps, in an appropriate
generalized sense.
 A version of the Douady-Hubbard theorem is valid in this situation: 
any cubic-like map
is quasi-conformally conjugate to a cubic polynomial with one escaping
critical point. It follows that all real cubic-like 
 Fibonacci maps are quasi-symmetrically conjugate.
 So, any example of a cubic-like Fibonacci map with unbounded geometry 
shows that all of them have unbounded geometry. Finally, we renormalize
a quadratic-like Fibonacci map into a cubic-like one which completes the proof
for the polynomial-like case.

In the last \S 8 we show that the limits of the renormalizations of a smooth
Fibonacci map are actually polynomial-like which completes the proof of
the Theorem.

{\bf Remark 4.} The Fibonacci recurrence is a well-known phenomenon
for monotone maps of the circle with golden rotation number.
The scaling laws in this situation were studied by Herman (at least
implicitly), by
Swiatek [Sw1] (smooth homeomorhisms with critical points) and by Tangerman and
Veerman [TV] (maps with flat spots). In the 
two former cases one has bounded geometry,
in the latter the geometry goes down to zero in the similar manner as in our
example.
Such circle maps are explicitely related to certain
unimodal maps of the interval which are different from ours but
also have a sort of Fibonacci recurrence ; see [PTT].
\smallskip

The notation $f^n$ will always be
used for the $n$-fold iterate of $f$.

Acknowledgement. We want to thank Branner, Douady, Sullivan,
and Tresser, for helpful conversations. We also profited  from the
discussions with the participants of the Stony Brook dynamical systems seminar,
particularly:
Brucks, Yu.Lyubich, Shishikura, Tangerman and Veerman.  
\bigskip\bigskip

\eject
\centerline{\bf \S2. Kneading.}\medskip

Let $f:I\to I$ be a unimodal map with minimum at $x=0$. As usual, let\break
$0=x_0\mapsto x_1\mapsto\cdots\;$ be the critical orbit, and
let $$ u(1)=1\,,\;\;u(2)=2\,,\;\;u(3)=3\,,\;\;u(4)=5\,,\;\;\ldots $$
be the Fibonacci numbers. In order to avoid the
hypothesis that $f$ is an even function, we will use the notation $x\mapsto x'$
for the order reversing involution, defined on some suitable subinterval
of $I$, which satisfies $f(x')=f(x)$. Let $\|x\|$ be the larger of $x$
and $x'$.\smallskip

{\bf Definition.} We will say that $f$ is a {\bit Fibonacci
map~} if $\|x_{u(n)}\|>\|x_{u(n+1)}\|$ for $n\ge 1$, so that
$$	\|x_1\|>\|x_2\|>\|x_3\|>\|x_5\|>\|x_8\|>\|x_{13}\|>\cdots \;,
\eqno (2-1) $$
and if $x_4<0\,$.

{\QP{\bf Lemma 2.1.} \it The map $f$ is a Fibonacci map if and only if
the signs of the successive images $x_i$ are given by
$$\eqalignno{	&\sgn(x_j)\;=\;\sgn(x_{j-u(n)})\qquad{\rm for}\qquad u(n)<j<u(n+1)\;,
	\quad {\rm with}&(2-2)\cr
	&\sgn(x_{u(n)})\;=\;(-1)^{(n+1)(n+2)/2}\;. &(2-3)}$$\medskip}

{\bf Remark 1.} Some condition such as $x_4<0$
is needed in order to avoid the\break uninteresting case
$$	x_1\;<\;0\;<\;\lim_{m\to \infty} x_m\;<\;\cdots\;<\;x_5\;<\;x_4\;<\;
 x_3\;<\;x_2\;. $$
(Note that such a map would have to have at least three fixed points,
counted with\break multiplicity. Thus
this particular case can never occur for a quadratic map.)\smallskip

{\bf Remark 2.} We can describe these conditions in different
language as follows. If we assume that $x_1<0<x_2\,$, then Conditions (2-2)
and (2-3) are completely equivalent to the statement that
the interval between $0$ and $x_{u(n)}$ is mapped homeomorphically by
the iterate $f^{\circ i}$ for $0\le i\le u(n-1)$, but
is not mapped homeomorphically by $f^{\circ u(n-1)+1}$.
The condition that some large iterate of $f$ restricted to
an interval $[a,b]$ is a homeomorphism is an invariant way of specifying
that $a$ is very close to $b$. Thus Lemma 2.1 can be thought of as
giving an invariant description of just how close $x_{u(n)}$ is to the
critical point.\smallskip

{\bf Remark 3.} The Branner-Hubbard description of $f$
would be rather different.
Following Yoccoz, they cut the interval not at the critical point, but
rather at the interior fixed point $\alpha<0\,$. In terms of the resulting
partition of the interval, the appropriate description of the critical orbit
is that
the two images $x_i$ and $x_{i+u(n)}$ lie on the same side of $\alpha$
for $i<u(n+1)-2\,$, but on opposite sides of $\alpha$
for $i=u(n+1)-2\,$.\bigskip

{\bf Proof of 2.1.} If (2-2) and (2-3) are satisfied, then according to
Remark 2 above, we see that the successive images $x_{u(n)}$ are closer and
closer to zero. Since $x_4<0$, it follows
that $f$ is a Fibonacci map. Conversely,
the proof that every Fibonacci map satisfies (2-2) and (2-3)
will be by induction on $n$, using the following induction hypothesis.

{\QP{\bf Hypothesis ${\bf H}_n$.} {\it For $i$ in the range $0<i<u(n)$ with
$ i\ne u(n-1)\,$, the  points $x_i$ have
sign as specified in Conditions $(2-2)$ and $(2-3)$
above, and furthermore $\|x_i\|>\|x_{u(n-1)}\|\,$.}\medskip}

\noindent The following elementary observation will be used over and over.
{\bit For any unimodal map with minimum at $x_0=0\,$:\hfil\break
\vskip -.1in
\indent\indent if $\qquad \|x_p\|<\|x_q\|\qquad$ then $\qquad x_{p+1}<x_{q+1}
\;\;. $}\hfil\break
\vskip -.1in
\noindent To start the induction, we must show that every Fibonacci map
satisfies ${\bf H}_4$. Since $\|x_1\|>\|x_2\|>\|x_3\|$ by definition,
we need only show that
$$x_1\,,\,x_4<0<x_2\;,\qquad{\rm and\;that}\qquad \|x_4\|>\|x_3\|\,.$$
Note first that the $\|x_i\|$ must
all be distinct. For otherwise the critical orbit would have only finitely
many distinct elements. We have assumed that $x_4<0\,$. 
If  $0<x_1$ then we see inductively that
$0<x_1<x_2<\cdots\,$, which contradicts our hypothesis. Similarly,
if $x_2<0$ hence $x_1<x_2<0$, then we see inductively that
$$x_1<x_3<x_5<\cdots<x_6<x_4<x_2<0\,,$$
which contradicts our hypothesis.
Finally, suppose that $\|x_4\|<\|x_3\|\,$. Applying the map $f$, we see
that $x_5<x_4<0$, and applying $f$ again we see that
$x_5<x_6\,$. Since
$\|x_5\|<\|x_3\|$ by hypothesis, hence $x_6<x_4\,$, we have $x_5<x_6<x_4<0\,$ and
a similar inductive argument shows that $x_5<x_7<x_9<
\cdots<x_8<x_6<x_4<0\,$, which again contradicts our hypothesis.
This proves ${\bf H}_4\,$.\smallskip

We will show that ${\bf H}_n\Rightarrow {\bf H}_{n+1}$ for
$n\ge 4\,$. Since $0<\|x_{u(n)}\| <\|x_{u(n-1)}\|$, we have
$$	x_1\;<\;x_{1+u(n)}\;<\;x_{1+u(n-1)}\;. $$
Now $x_i$ and $x_{i+u(n-1)}$ have the same sign for $0<i<u(n-2)$ by
${\bf H}_n$. Hence it follows by induction on $i$ that $x_{i+u(n)}$ lies
between them, and hence also has the same sign, for $i$ in this range.
Since both $x_i$ and $x_{i+u(n-1)}$ have absolute value greater than
$\|x_{u(n-1)}\|$ by ${\bf H}_n\,$, it follows also that
$\|x_{i+u(n)}\|>\|x_{u(n-1)}\|>\|x_{u(n)}\|\,$, for $i$ in this range.
For $i=u(n-2)$, this argument
proves that $x_{u(n-2)+u(n)}$ lies between $x_{u(n-2)}$ and $x_{u(n)}$, but
does not determine its sign. However, it does follows that
$$ 0<\|x_{u(n-2)+u(n)}\|<\|x_{u(n-2)}\|\;,\qquad{\rm hence}\qquad x_1<
	x_{1+u(n-2)+u(n)}<x_{1+u(n-2)}\;. $$
Now a similar inductive argument shows that $x_{i+u(n-2)+u(n)}$ lies between
$x_i$ and\break
$x_{i+u(n-2)}\,$, and hence has the required sign, for $0<i<u(n-3)$.
Furthermore, this shows that $\|x_{i+u(n-2)+u(n)}\|>\|x_{u(n-1)}\|>\|x_{u(n)}\|$
for $i$ in this range.
In the limiting\break case $i=u(n-3)$, this argument proves that $x_{i+u(n-2)
+u(n)}=x_{u(n+1)}$ lies between\break $x_{u(n-3)}$ and
$x_{i+u(n-2)}=x_{u(n-1)}$,
but does not determine its sign. However, since\break
$\|x_{u(n+1)}\|<\|x_{u(n-1)}\|
<\|x_{u(n-3)}\|\,$, this proves that $x_{u(n-3)}$ and $x_{u(n-1)}$ have opposite
sign, so that $x_{u(n-1)}$ also has the required sign.
Thus, we have almost proved ${\bf H}_{n+1}$. The only missing pieces
of information are the sign and magnitude of $x_i$ for $i=u(n-2)+u(n)\,$.
\smallskip

We must prove that $\|x_{u(n-2)+u(n)}\|>\|x_{u(n)}\|\,$. But if
$\|x_{u(n-2)+u(n)}\|<\|x_{u(n)}\|$ then
$$	x_1\;<\; x_{1+u(n-2)+u(n)}\;<\;x_{1+u(n)}\;. $$
This is impossible. For a similar inductive argument would show that
$x_{i+u(n-2)+u(n)}$ must be between $x_i$ and $x_{i+u(n)}$ for $0<i\le u(n-2)$.
In particular, taking
$i=u(n-3)$ it would follow that $x_{u(n+1)}$ must be
between $x_{u(n-3)}$ and $x_{u(n-3)+u(n)}$. By the part of ${\bf H}_{n+1}$
which has already been proved, these two have the same sign, and it would
follow that $\|x_{u(n+1)}\|>\|x_{u(n)}\|\,$, which contradicts our hypothesis.
Thus $\|x_{u(n-2)+u(n)}\|>\|x_{u(n)}\|\,$.
\smallskip

Now recall that $x_{u(n-2)+u(n)}$ is known to lie between $x_{u(n-2)}$ and
$x_{u(n)}$. Since $\|x_{u(n-2)+u(n)}\|>\|x_{u(n)}\|\,$, it follows easily
that $x_{u(n-2)+u(n)}$ has the same sign as $x_{u(n-2)}$. This
completes the proof that ${\bf H}_n\Rightarrow {\bf H}_{n+1}$. \QED
\medskip

To show that this result is not vacuous, we must prove the following.

{\QP{\bf Lemma 2.2.} \it Fibonacci maps exist.\medskip}

We will outline two different proofs. The proof below is an
immediate application of the formal machinery
of kneading theory, as developed in [MT]. 
An alternative proof, which is more direct and gives a more explicit
description of the critical orbit, will be given in Lemma 3.1. Both proofs
will make use of the following.\smallskip

{\bf Definition 2.3.} By a {\bit Fibonacci sum~} we will mean a finite or
infinite formal sum $$ \mu\;=\;u(n_1)+u(n_2)+u(n_3)+\cdots $$
of {\bit non-consecutive\/}
Fibonacci numbers. That is, we always assume that $n_{i+1}\ge n_i+2\,$,
with $n_1\ge 1\,$. It is not difficult to check that every positive integer has
a unique expression as a finite Fibonacci sum. As an example, the difference
$u(n)-1$ can be expressed as
$$	u(n)-1\;=\;\left\{ \eqalign{
	u(1)+u(3)+u(5)+\cdots+u(n-1)&\qquad{\rm for}\;\; n\;\;{\rm even}\,,\cr
	u(2)+u(4)+u(6)+\cdots+u(n-1)&\qquad{\rm for}\;\; n\;\;{\rm odd}\,.
	}\right.\eqno (2-4) $$
(For infinite Fibonacci sums, compare the proof of Lemma 3.2.)
\smallskip

As in [MT], we describe the kneading invariant of a unimodal map $f$
by a formal power series $\;D(t)=1+\epsilon_1t+\epsilon_2t^2+\cdots\,$,
where each coefficient $\epsilon_n$ is equal to $+1$ or $-1$ according as
the function $x\mapsto |f^{\circ n}(x)|$ has a local minimum or local
maximum at the origin. Since the $x_i$ are non-zero for $i>0$, we can
check inductively that
$$	\epsilon_n\;=\;\sgn(x_1x_2\cdots x_n)\,. \eqno (2-5) $$
Such a kneading invariant is {\bit admissible\/} (ie., actually occurs)
if and only if the inequality
$$	\sum_0^\infty \epsilon_i\,t^i\;\le\;\sum_0^\infty\, (\epsilon_m
	\epsilon_{m+i})\,t^i	\eqno (2-6) $$
is satisfied for every $m\ge 1$. Here, by definition, an inequality
$\sum a_it^i <\sum b_it^i$ between formal power series means
that the first difference $b_i-a_i$ which is non-zero is
actually positive. Thus, for each $m$ we require that the smallest $i$
for which $\;\epsilon_{m+i}\;\ne\;\epsilon_m\,\epsilon_i\;$\break (if any
such exist) must satisfy $\epsilon_i=-1\,$.\smallskip

In the case of a Fibonacci map, it follows inductively from (2-2), (2-3) 
and (2-5)
that we must have $\epsilon_{u(n)}=-1$ for every Fibonacci number $u(n)$.
In fact, according to (2-5),
$\epsilon_{u(n+1)}$ is equal to $\epsilon_{u(n)}$ multiplied by the
sign of the product $x_{u(n)+1}x_{u(n)+2}\cdots x_{u(n+1)}$. This coincides
with $\sgn\big(x_1x_2\cdots x_{u(n-1)}\big)=\epsilon_{u(n-1)}=-1$ except that
the very last factor $x_{u(n-1)}$ has the wrong sign. Thus it follows
inductively that $\epsilon_{u(n)}=\epsilon_{u(n+1)}=-1$ for all $n$.
In other words, each map $x\mapsto |f^{\circ u(n)}(x)|$ must have a local
maximum at $x=0$. For a $k$-fold Fibonacci sum
$$	m\;=\;u(n_1)+\cdots+u(n_k)\,,\qquad{\rm where\;always}\quad
 n_1\ge 1\;\;{\rm and}\;\;n_{i+1}\ge n_i+2\,, \eqno (2-7) $$
Equations (2-2) and (2-5) imply that
$ \epsilon_m$ is equal to the product $\epsilon_{u(n_1)
+\cdots+u(n_{k-1})}\epsilon_{u(n_k)}\,. $ Hence it follows inductively
that $\epsilon_m=(-1)^k\,$. Thus, in order to prove 2.2
we need only show that the formal power series $\sum\epsilon_mt^m$, with
$\epsilon_m$ defined by this equation,
satisfies Condition (2-6). That is, for each fixed $m$ the smallest $i$ with
$\epsilon_{m+i}\ne\epsilon_m\epsilon_i$ must satisfy $\epsilon_i=-1\,$.
However, if we express $m$ as a Fibonacci sum as above,
then it is not hard to show that the smallest $i$ with $\epsilon_{m+i}\ne
\epsilon_m\epsilon_i$ is either $i=u(n_1-1)$ or $i=u(n_1)$ or (in the special
case $n_1=1$) $i=2$. Since $\epsilon_i=-1$ in each of these cases, the required
inequality (2-6) follows. This completes the proof of 2.2. \QED\smallskip

{\bf Proof of Theorem 1.1.} Since any unimodal kneading invariant which
is admissible can be
realized by a quadratic map, we can certainly find at least one quadratic
map $f_c$ which realizes the given kneading
invariant. (See for example [MT].) But for any real
quadratic map $f_c$ which is not infinitely renormalizable and has
no attracting periodic orbit, Yoccoz has recently shown that the constant
$c$ is uniquely determined by its kneading invariant. (This is an immediate
corollary of his much more general
result about complex quadratic parameter space.) Since
it is easy to check that a quadratic Fibonacci map is not renormalizable
and has no attracting periodic orbit, this proves 1.1. \QED
\bigskip\bigskip

\centerline{\bf \S3. The critical orbit.}\medskip

Out of the kneading data, it is not difficult to determine the precise ordering
of the points $x_m$ in the critical orbit. We can describe the resulting
ordering by a fairly concrete model as follows. The construction will provide an alternative proof of 2.2.\smallskip

Choose a parameter $0<t<1-t^2\,$, or in other words
$$0\;<\;t\;<\;(\sqrt 5-1)/2\;=\;.61803\cdots\,, $$
for example $t={1\over 2}\,$. Now for each integer
$m\ge 0$,
expressed as a Fibonacci sum (2-7), define a real number $y_m$ by the formula
$$ y_m\;=\;\pm\big(t^{n_1}-t^{n_2}+-\cdots\pm t^{n_k}\big)\;, $$
where the initial sign is to be $-1$ for $n_1\equiv 0\,,\,1\;(\mod 4)\,$,
and $+1$ for $n_1\equiv 2\,,\,3\;(\mod 4)\,$, as in (2-3) above. Thus the
initial term $\pm t^{n_1}$ is the dominant one, and subsequent terms
alternate in sign, decreasing by a factor of $t^2$ or more at each step since
$n_{i+1}\ge n_i+2\,$.\smallskip

{\bf Remark 3.1.} More precisely, this ordering can be described as
follows. For Fibonacci sums $m$ with different dominant
terms, the order of the $y_m$ is determined by the rules:
$$ y_{1+\cdots}<y_{5+\cdots}<y_{8+\cdots}<y_{34+\cdots}
 <\;\cdots\;<0	<\;\cdots\;<y_{21+\cdots}<y_{13+\cdots}
	<y_{3+\cdots}<y_{2+\cdots}\;\;. $$
Here, in each case, the dots in the subscript stand for higher
terms, which may be zero, for an arbitrary Fibonacci sum. For two
Fibonacci sums which have the same leading summands $u(n_1)+\cdots+u(n_k)$ but
differ at the $(k+1)$-st summand, the relative order is
determined as follows. Setting $s=u(n_1)+\cdots+u(n_k)$, we have
$$	|y_{s}|\;>\; \cdots\;>\; |y_{s+u(n_k+5)+\cdots}|\,>\,
	|y_{s+u(n_k+4)+\cdots}|\,>\,
	|y_{s+u(n_k+3)+\cdots}|\,>\, |y_{s+u(n_k+2)+\cdots}| $$
if $k$ is odd; and the same but with all inequalities
reversed if $k\ge 2$ is even. Here all of these
points $y_{s+\cdots}$ have the same sign, depending only on the leading
summand $n_1$,
as described above.\smallskip


We claim that the resulting ordering of the
$y_m$ is precisely the required ordering of the points $x_m$ in the critical
orbit. More precisely, we will prove the following.

{\QP{\bf Lemma 3.2.} \it The correspondence $y_m\mapsto y_{m+1}$
is unimodal, that is, it is\break monotone increasing on the set of $y_m$
for which $y_m\ge 0$, but monotone\break
decreasing for $y_m\le 0$. Furthermore,
this correspondence is uniformly continuous. Thus, if we extend
linearly over each gap between the $y_m$, then we obtain a continuous
unimodal map $F$ from the interval $[y_1\,,\,y_2]$ to itself, satisfying the
Fibonacci condition that
$$	y_1<y_2'<y_3'<y_5<y_8<y_{13}'<\cdots<0 \;, $$
where $y_m=F^m(0)$.
$($Here, as in \S2, we use the notation $y\mapsto y'$
for the orientation reversing involution of the subinterval $[y_2'\,,\,y_2]$
which satisfies the condition that $F(y')=F(y)\;.)$\medskip}

{\bf Proof.}
It is convenient to divide the various $y_m$ into intervals $A_n\,,\; n\ge
 0\,$, which are ordered according to the following pattern:
$$ A_2<A_6<A_{10}<\cdots<A_8<A_4<A_0\;\;\le\;\; A_1<A_5<A_9<\cdots<A_{11}
	<A_7<A_3\;.$$
(Here the two sequences $\{A_{2n}\}$ and $\{A_{2n+1}\}$ converge towards the
two pre-images of zero. Compare 
 3.6.)
Let $A_0=[y_5\,,\,0]$ be the closed interval containing
all $y_{u(n)+\cdots}$ with\break $n\equiv 0\,,\, 1\;(\mod 4)\,,\; n\ge 4\,$,
and also containing the limit
point zero. Here, as above, the notation $u(n)+\cdots$ stands for an
arbitrary Fibonacci sum with leading term $u(n)$.
Similarly, let $A_1=[0\,,\,y_3]$ be the interval containing all
$y_{u(n)+\cdots}$ with $n\equiv 2\,,\, 3\;(\mod 4)\,,\; n\ge 3$, together
with the limit point zero. For $n\ge 2$ even, let $A_n$ be the smallest
interval\break
containing all $y_m$ with $m$ of the form $u(1)+u(3)+\cdots+u(n-1)+{\rm (higher
\;terms})$, where the higher terms if any must start with $u(n+2)$ or higher.
Using the identity (2-4),
it follows easily that $A_n$ is equal to the closed interval spanned by the
two points
$y_{u(n)-1}$ and $y_{u(n)+u(n+2)-1}$. Here the relative order of these
two endpoints depends on whether $n$ is congruent to 0 or 2 modulo 4.
Similarly, for $n\ge 3$ odd, we define $A_n$ to be the smallest
interval containing all $y_m$ with $m$ of the form $u(2)+u(4)+\cdots+u(n-1)
+({\rm higher})$,
where again the higher summands if any must start with $u(n+2)$ or higher.
Again using the identity (2-4),
we see that this interval $A_n$ is again spanned by the points $y_{u(n)-1}$
and $y_{u(n)+u(n+2)-1}$, where the relative order of the two end points
depends on whether $n$ is congruent to 1 or 3 modulo 4.

It is not difficult to show that
every $y_m$ with $m>0$ belongs to exactly one of these intervals,
and that these points are ordered according to the pattern described above.
For $y_m\in A_n$ a brief computation shows that the map $y_m\mapsto
 y_{m+1}$ is linear with slope $(-1)^{n-1}$. In particular, it is either
order preserving or order reversing according as\break $A_n\subset[0,y_2]$ or
$A_n\subset[y_1\,,\,0]$. If we extend this map to be
linear in the gap between $A_n$ and $A_{n+4}$, then computation shows that
the slope in this gap takes the value
$$	{\Delta F(x)\over\Delta x}\;=\; (-1)^{n-1}{t^n-t^{n+2}-t^{n+4}
	\over t^{n+1}-t^{n+2}-t^{n+3}} $$
for $n>0\,$. This is independent of $n$ except for sign.
For $n=0$ it takes a different value, but still with the appropriate 
(negative) sign. As an example, for
$t={1\over 2}$ this gap slope is equal to $\pm{11\over 2}$ for $n>0$, and is
$-{6\over 5}\,$ for $n=0\,$. In this way, we obtain the required explicit
unimodal map $F$ which realizes the given kneading data. This 3.2,
and completes the alternate proof of 2.2. \QED\smallskip

{\QP{\bf Lemma 3.2.} \it If the Fibonacci map $f$ has no ``homtervals''
within the interval $[x_1\,,\,x_2]$,
that is, if the pre-critical points are everywhere dense, then
$f$ restricted to this interval is topologically conjugate
to this model map $F$.\medskip}

\noindent The proof is straightforward. \QED

{\bf Remark.} By definition, a {\bit homterval\/} is a subinterval of
$I$ which is mapped homeomorphically by all iterates of $f$. A
{\bit wandering interval\/} is a homterval which is not contained in
the basin of attraction for any periodic orbit.
According to Guckenheimer [G1], a unimodal map has no wandering intervals
within $[x_1\,,\,x_2]$ provided that it has
negative Schwarzian, with non-flat
critical point. According to de Melo and van Strien [MS], it has no wandering
intervals provided that it is
sufficiently smooth, with non-flat critical point.
(See also Blokh and Lyubich [L], [BL1].)
$\qquad\qquad$\medskip

{\QP{\bf Lemma 3.4.} \it More generally, if a Fibonacci map has no wandering
intervals, then its critical orbit closure $\overline{\cal O}$ is a
Cantor set, homeomorphic to the corresponding critical orbit closure for
the model map $F$. In particular, this Cantor set is canonically
homeomorphic to the set of all finite or
infinite Fibonacci sums, suitably topologized.\medskip}

{\bf Proof of 3.4.}
The appropriate topology for the set of all finite or
infinite\break Fibonacci sums can be described as follows. Let $\Sigma$ be
the ``Fibonacci shift'', consisting of all sequences $(a_1\,,\,a_2\,,\,\ldots)$
of zeros and ones with no two consecutive ones. 
In other words, $\Sigma$ is a one-sided subshift of finite type corresponding
to the matrix $T=\pmatrix{1&1\cr1&0}$.
(The name is suggested
since the number of cylinders   in $\Sigma$ of length $n$
is equal to $u(n+1)$.) This set $\Sigma$ is topologized as a
subset of the infinite Cartesian product $\{0,1\}\times\{0,1\}\times\cdots$.
Each sequence $\{a_n\}\in\Sigma$ determines an associated 
Fibonacci sum  $\mu=\sum a_nu(n)$, and we
give the set consisting of all Fibonacci sums the corresponding compact
topology. It is
easy to check that the correspondence $m\mapsto x_m\,$, where $m$
ranges over positive
integers expressed as finite Fibonacci sums, extends uniquely to a
homeomorphism $\mu\mapsto x_\mu\in\overline{\cal O}$, where now $\mu$
ranges over finite or infinite Fibonacci sums. Further
details of the proof are straightforward. \QED

{\bf Remark 3.5.} It is sometimes convenient to
partially order the Cantor set $\Sigma$ using
lexicographical order from the right. Thus two sequences of zeros and ones,
with no two consecutive ones, are comparable whenever they are eventually
equal, or in other words have the same
tail. In terms of this ordering,
the map from $\Sigma$ to itself which corresponds to the map
$f|\overline{\cal O}$ can be described as the {\bit immediate successor
function\/}, which carries each such sequence to the next largest sequence with
the same tail
(such a transformation is called an {\it adic shift}, compare [V]). 
However, there are two exceptional sequences which are maximal,
and hence have no successor, immediate or otherwise, namely the two
sequences $(1,0,1,0,\ldots\,)$ and $(0,1,0,1,\ldots\,)$ corresponding to
the Fibonacci sums $\;1+3+8+\cdots\;$ and $\;2+5+13+\cdots\;$ respectively.
These both map to the zero sequence. (Compare (2-4) in \S2.)

{\QP{\bf Corollary 3.6.} \it
The mapping $f$ from the Cantor set $\overline{\cal O}$
onto itself is one-to-one except that the point zero has two different
pre-images, corresponding to the infinite Fibonacci sums $u(1)+u(3)+u(5)+
\cdots$ and $u(2)+u(4)+u(6)+\cdots\,$.\medskip}

\noindent The proof is straightforward.\QED\smallskip

Here is a more explicit description of this Cantor set as a subset of the
real line. For each $n\ge 1$ let $I^n\subset\R$ be the smallest closed
interval containing all of the points $x_{u(q)}$ with $q\ge n$. Thus $I^n$
is a closed neighborhood of the origin. One
end point of this interval is $x_{u(n)}$ and the other end point is either
$x_{u(n+1)}$ or $x_{u(n+2)}$ according as $n$ is odd or even. Note that the
map $f$ folds $I^n$ over onto the closed interval $[x_1\,,\,x_{u(n)+1}]$,
which in turn maps onto the closed interval $[x_{u(n)+2}\,,\, x_2]$
provided that $n\ge 3$.
For each $k\ge 0$, we will use the notation $I^n_k$ for the image
$f^k(I^n)$. According to \S2, this image $I^n_k$ is disjoint
from the origin for $1\le k<u(n-1)$, but contains the origin for $k=u(n-1)$.
However, $I^n_{u(n-1)}$ contains a smaller interval $I^{n+1}_{u(n-1)}$
which again is disjoint from the origin.
It will be convenient to use the notation
$$	J^n=I^{n+1}_{u(n-1)}\;,\quad{\rm and\;more\;generally}\quad
 J^n_k=f^k(J^n)
	=I^{n+1}_{k+u(n-1)}\; .$$
Note in particular that $J^n_{u(n-2)}=I^{n+1}_{u(n)}$.\smallskip

{\bf Definition.} Let $M^n$ be the $u(n)$-fold union
$$	M^n\;=\;\bigcup_{0\le k<u(n-1)} I^n_k\;\;\;\cup\;
	\bigcup_{0\le k<u(n-2)} J^n_k\;. $$
For example (listing the subintervals from left to right):
$$\eqalign{ M^1\;=\;& [x_1\,,\,x_2]\cr
	M^2\;=\;& [x_1\,,\,x_4]\cup[x_5\,,\,x_2]\cr
	M^3\;=\;& [x_1\,,\,x_4]\cup[x_5\,,\,x_3]\cup[x_7\,,\,x_2]\cr
	M^4\;=\;& [x_1\,,\,x_6]\cup[x_{12}\,,\,x_4]\cup[x_5\,,\,x_{13}]\cup
	[x_{11}\,,\,x_3]\cup[x_7\,,\,x_2]\cr} $$
and so on. 

{\QP{\bf Lemma 3.7.} \it The $u(n)$ closed intervals 
$$	 I^n_0\;,\; I^n_1\;,\;\ldots\;,\;
	 I^n_{u(n-1)-1}\qquad{\rm and}\qquad
	J^n_0\;,\;\ldots\;,\; J^n_{u(n-2)-1}\; $$
are pairwise disjoint. Denoting their union by $M^n$ as above, the
$M^n$ form a nested sequence of
closed sets $M^1\supset M^2\supset M^3\supset\cdots$  with
intersection equal to the Cantor set $\overline{\cal O}$.
\medskip}

{\bf Proof}. We will show by induction on $n$ that the $u(n)$ subintervals
of $M^n$ are pairwise disjoint, that the $M^n$ are nested, and that
each $M^n$ contains the critical orbit closure. The idea of the proof is to
show that, as we pass from $M^n$ to $M^{n+1}$, each of the $u(n-1)$ intervals
$I^n_k\subset M^n$ will be replaced by two subintervals $I^{n+1}_k$ and
$J^{n+1}_k$ in $M^{n+1}$, while each of the $u(n-2)$ intervals $J^n_k
=I^{n+1}_{k+u(n-1)}$ remains unchanged.\smallskip

To start the induction, it is trivially true that
$M^1=[x_1\,,\,x_2]$ contains the critical orbit closure. The first
step in the induction is to note that each $I^n$ contains $I^{n+1}$
and $J^{n+1}$ as disjoint subsets. For example if $n\equiv 3\;(\mod 4)$
then these two subinterval of $$I^n=[x_{u(n+1)}\,,\,x_{u(n)}]$$
are situated as follows:

\centerline{\psfig{figure=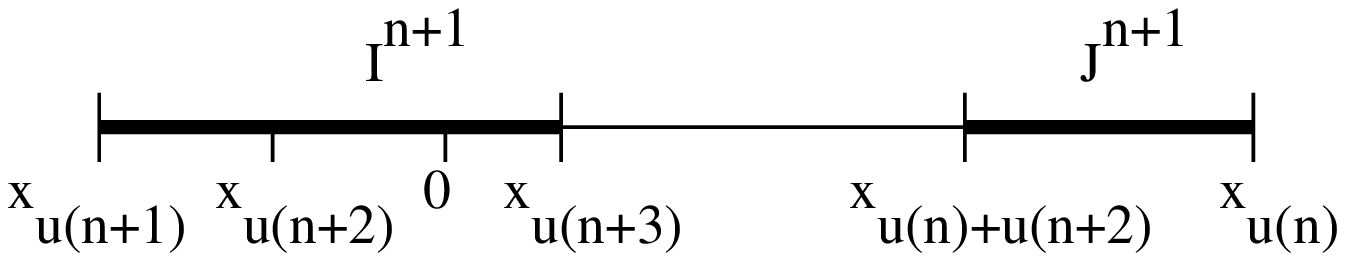,height=1.5in}}
\centerline{Figure 1. The interval $I^n$ in the case $n\equiv 3\;(\mod 4)$.}
\bigskip

\noindent The picture for $n\equiv 1\;(\mod 4)$ is a mirror image, and the
pictures for $n\equiv 0\,,\,2\;(\mod 4)$ are quite
similar. Note that the map $f^{ u(n)}$ folds the subinterval $I^{n+2}
\subset I^{n+1}$
over onto $J^{n+1}$, while the map $f^{u(n-1)}$ carries $J^{n+1}$ back onto
a neighborhood of the\break origin, spanned by the two points $x_{u(n+1)}$ and
$x_{u(n+3)}$. In the case $n\equiv 3\;(\mod 4)$\break as illustrated, $I^{n+2}$
is the interval $\;[x_{u(n+2)}\,,\, x_{u(n+3)}]\,$, while the image
$\;f^{\circ u(n-1)}(J^{n+1})\break =\;[x_{u(n+1)}\,,\,x_{u(n+3)}]\;$
coincides with the interval $\;I^{n+1}$.\smallskip

It follows easily from 3.1 and 3.4 that the two subintervals
$$	I^{n+1}\,,\;J^{n+1}\;\subset\;I^n $$
are indeed disjoint, and together contain all of the points of $\overline
{\cal O}\cap I^n$. For $1\le k<u(n-1)$, a similar
argument shows that the two subintervals
$$	I^{n+1}_k\,,\;J^{n+1}_k\;\subset\;I^n_k $$
are disjoint, and together contain all of the points of $\overline
{\cal O}\cap I^n_k$. This completes the induction, and shows that
$$	M^1\supset M^2\supset M^3\supset\cdots\supset\overline{\cal O}\;.$$
Since each endpoint of each subinterval of $M^n$ belongs to the orbit
${\cal O}$, using the hypothesis that there are no wandering intervals we
see easily that $\;\bigcap M^n\;$ is equal to $\;\overline{\cal O}$. \QED

Using the sets $M^n$ one can give another description of the above correspondence
between $\overline{\cal O}$ and $\Sigma$ (see 3.4). 
Given $x\in \overline{\cal O}$,
let $M^n(x)$ be an interval of the set $M^n$ containing $x$. Then set $a_n=0$
if $M^n(x)=I^n_k$, and $a_n=1$ if $M^n(x)=I^{n+1}_k$ (for appropriate $k$'s).
One can check that $\{a_n\}\in \Sigma$ is the sequence corresponding to 
$x\in \overline{\cal O}$.

In what follows we will use the notation $M^n_{a_1...a_n}$ for the interval
of $M^n$ corresponding to the cylinder $[a_1...a_n]\subset\Sigma$.

{\QP{\bf Lemma 3.8.} \it Still assuming that there are no
wandering intervals, the points $x_i\;,\;\;i\ge 1$ are
the endpoints of the complementary intervals for the critical\break
orbit closure
$\overline{\cal O}\subset\R$. More explicitly, the Cantor set
$\overline{\cal O}$ can be obtained from the closed interval $[x_1\,,\,x_2]$
by removing a dense collection of disjoint open sub\-intervals $(x_p\,,\,x_q)$
as follows. If one of $p\,,\,q$ is a Fibonacci sum of the form
$$u(n_1)\;+\;\cdots\;+\;u(n_{k-1})\;+\;u(n_k)\;+\;u(n_k+2)$$
with $k\ge 2$, then the other is equal to
$$u(n_1)\;+\;\cdots\;+\;u(n_{k-1})\;+\;u(n_k+1)\,.$$
On the other hand, if one is $\;u(n)
+u(n+2)\;$, then the other is either\break
$u(n+1)\;$ or $\;u(n+3)\;$ according as $n$
is even or odd.\medskip}

\noindent As an example, the first seven open subintervals to be removed are as
follows, in their natural order:
$$	(x_9\,,\,x_{19})\,\cup\,(x_6\,,\,x_{12})\,\cup\,(x_4\,,\,x_5)
\,\cup\,(x_{18}\,,\,x_8)\,\cup\,(x_{13}\,,\,x_{11})\,\cup\,
(x_3\,,\,x_7)\,\cup\,(x_{20}\,,\,x_{10})\;. $$
In other words, the Cantor set $\overline{\cal O}$
is contained in the disjoint union
$$ [x_1\,,\,x_9]\,\cup\,[x_{19}\,,\,x_6]\,\cup\,[x_{12}\,,\,x_4]\,\cup\,
[x_5\,,\,x_{18}]\,\cup\,[x_8\,,\,x_{13}]\,\cup
\,[x_{11}\,,\,x_3]\,\cup\,[x_7\,,\,x_{20}]\,\cup\,[x_{10}\,,\,x_2] $$
(which coincides with the closed set $M^5$).
The proof of this statement is a straightforward consequence of the ordering
of the points in the critical orbit, as described above. \QED
\smallskip

We can obtain a different model for this critical orbit closure as follows.
Let $$\gamma\;=\;(1-\sqrt 5)/2\;=\;-.61803\cdots\,,$$
so that $\gamma=\gamma^2-1$. To each finite
or infinite Fibonacci sum $\mu=u(n_1)+u(n_2)+\cdots$,\break
let us assign the real
number $\;\phi(x_\mu)\,=\,\gamma\,(\gamma^{n_1}+\gamma^{n_2}+\cdots\;)\;$
modulo one.

{\QP{\bf Lemma 3.9.} \it The resulting map $\phi$ from the critical orbit
closure $\overline{\cal O}$ onto
the circle $\R/\Z$ is one-to-one except at the countably many iterated
pre-images of zero. It semi-conjugates the map $f|\overline{\cal O}$ onto
the golden rotation $\;t\mapsto t+\gamma\;(\mod 1)\;$.\medskip}

{\bf Proof.} It is easy to check that $\phi$ is well defined and continuous.
Note that the identity $u(n-1)+u(n)=u(n+1)$ corresponds to the identity
$\gamma^{n-1}+\gamma^n=\gamma^{n+1}$. Using this fact, it is not difficult
to check the required identity
$$	\phi(f(x_\mu))\;=\;\phi(x_{\mu+1})\;\equiv\;
	\phi(x_\mu)+\gamma\quad(\mod 1)\,.  $$
Thus the image is a compact subset of the circle, invariant under the
golden rotation, and hence is equal to the entire circle.
Now consider any Fibonacci sum with leading term $u(n)$. A brief computation
shows that the corresponding image
$$	\phi(x_{u(n)+\cdots})\;=\;\gamma^{n+1}+\cdots $$
lies somewhere between
$$	\gamma^{n+1}+\gamma^{n+3}+\gamma^{n+5}+\cdots\;=
\; \gamma^{n+1}/(1-\gamma^2)\;=\;-\gamma^n$$
and
$$	\gamma^{n+1}+\gamma^{n+4}+\gamma^{n+6}+\cdots\;=\;\gamma^{n+1}
-\gamma^{n+3}\;=\;-\gamma^{n+2}\;. $$
Thus, depending on the leading summand, the image $\phi(x_\mu)$ lies in
one of the non-overlapping intervals
$$ [-\gamma^2,-\gamma^4]\cup[-\gamma^4,-\gamma^6]\cup[-\gamma^6,-\gamma^8]\cup
 \cdots	\cup\{0\}\cup\cdots\cup[-\gamma^7,-\gamma^5]\cup[-\gamma^5,-\gamma^3]
	\cup[-\gamma^3,-\gamma]\;, $$
having total length $-\gamma-(-\gamma^2)=1\,$.
Hence the value $\phi(x_\mu)\in\R/\Z$ determines the leading summand
$u(n)$ uniquely, except in countably many cases which can be explicitly
described. For two Fibonacci sums with the same leading
term, a similar argument shows that
the value $\phi(x_\mu)$ determines the second term uniquely, again with the
exception of countably many cases which can be explicitly described;
and a similar argument applies to higher terms. \QED\medskip

{\QP{\bf Corollary 3.10.} \it With hypotheses as above, the map $f|
\overline{\cal O}$ is minimal, that is every orbit is dense, and has
topological entropy zero. Furthermore this map is uniquely ergodic,
that is it has one and only one invariant probabiltity measure.\medskip}

{\bf Proof.} This follows easily from the corresponding assertion for
an irrational rotation of the circle. \QED\smallskip

Combining 3.4--3.10, this evidently completes the proof of Theorem 1.2. \QED
\smallskip

\bigskip\bigskip\bigskip

\centerline{\bf \S4. A priori bounds.}\bigskip

In the following two sections we assume that $f:[-1,1]\rightarrow [-1,1]$
 is a $C^2$-smooth unimodal map 
with non-degenerate minimum point 0, and normalized by the condition
$f(-1)=f(1)=1$ (which does not restrict the generality). 
 Denote this class of maps by $\cal U$,
and let us discuss topology on this space.

We will mainly be interested in the subspace ${\cal U}_0\subset{\cal U}$
consisting of those $f$ for
which $f$ is an even function, $f(-x)=f(x)$. We will first discuss the
differentiability conditions and topology on this subspace, and then
generalize to the full space ${\cal U}$. If $f$ is even, then we can
write it uniquely as $$f(x)=A_{x_1}\circ g\circ Q$$ where  $Q$
is the squaring map
$\xi\mapsto\xi^2$, $g$ is some orientation preserving diffeomorphism
of $[0,1]$, and $A_{x_1}$ is the orientation preserving affine map which carries
$[0,1]$ onto $[x_1,1]$, where $x_1=f(0)$ is the critical value. 

Now  the $C^k$-topology on ${\cal U}_0$, $k\leq 2$,
 comes from the $C^k$-topology on the space
of diffeomorphisms $g$ , together with the line topology on
the range  of the parameter $ x_1$. 
Let $\| f\|$ denote the maximum of the $C^2-$norms for
$g , g^{-1}$ which is a continuous functional in $C^2$-topology on our
 space.

To obtain a corresponding topology of the full space $\cal U$ we need one extra
step. Let $x\mapsto x'$ be the orientation reversing diffeomorphism of
$T$ which satisfies $f(x)=f(x')$.  This involution is certainly $C^2$-smooth.
 Consider a map  $B: x\mapsto (x-x')/2$.
 Evidently $f$ can be expressed as a function of
$(x-x')^2/4$, so that we have a presentation
$f(x)=A_{x_1}\circ g\circ Q\circ B$ instead of the above one.
Now we must
incorporate the $C^k$ topology on the involution as part of our topology.
In practice, it is easiest simply to carry out this symmetrizing change
of coordinate $x\mapsto (x-x')/2$
in the beginning, and thereafter to deal
only with even maps $f$. Moreover, we can also assume without loss of
generality that $f$ is purely quadratic $x\mapsto x^2-c$ near 0
(since any $f\in {\cal U}$ is $C^2$-conjugate to such one).

Denote by ${\cal F}$ the subspace of Fibonacci maps $f\in {\cal U}$.

The following notations will be kept throughout the paper:

$$d_n=|x_u(n)|\; ,\qquad \lambda_n=d_n/d_{n-1}\;.$$
The goal of this section is to obtain some a priori estimates for
the $\lambda_n$ (compare [G2], [L], [MMSS], [BL3], [M],
 [S],...). The proofs are based upon the Schwarz lemma
and the Koebe Principle stated in the Appendix.

First let us introduce a convenient terminology and notations. A family of 
intervals ${\bf G}=\{G_i\}_{i=0}^n$ is called a {\it chain of intervals}
if $G_i$ is a component of $f^{-1} G_{i+1}$ for $i=0,1,...,n-1$. The chain is
called {\it monotone } if all maps $f: G_i\rightarrow G_{i+1}$ are homeomorphisms.

For a given interval $G$ and a point $x$ such that $f^n x\in G$ one can 
construct a chain $G_0 , G_1 ,..., G_n\equiv G$ {\it pulling G back} along
the $n$-orbit of $x$. This construction is an efficient tool in one
dimensional dynamics because it is often possible to estimate the distortion
of $f^n$ along chains of intervals (see [L] , [S] ).

For a family of intervals ${\bf G}=\{G_i\}$ denote by $|{\bf G}|=\sum |G_i|$
 the total length of intervals
$G_i$ and by mult${\bf G}$ the maximal intersection multiplicity of intervals
$G_i$, that is the maximum number of $G_i$  having non-vacuous
intersection.

Let us consider now the pull-back
 $${\bf H}^{n+1}=\{ H^{n+1}_m\}_{m=0}^{u(n)-1},
\qquad H^{n+1}_0\equiv H^{n+1}\supset I_1^{n+1} \eqno  (4-1)$$
of the interval $T^{n-2}$ 
along the orbit $\{f^mI_1^{n+1}\}_{m=0}^{u(n)-1}$.
 The following two topological
lemmas easily follow from the above combinatorics.

{\QP{\bf Lemma 4.1.} The chain ${\bf H}^{n+1}$ is monotone
(so that $f^n$ monotonously maps $H^{n+1}$ onto $T^{n-2}$). \medskip}  

Let us consider any interval $I=I_k^l\; , l\in\{n, n+1\}$, of the family $M^n$
different from $I_0^n, I_1^n, I_2^n$. Define an interval 
$F\equiv F_n(I)\supset I$ as follows

(i) If $I\neq J^n$ then $F$ is the convex hull of two neighbors of $I$ in the
family $M^n$;

(ii) If $I=J^n$ then $F$ is the half of the interval $T^{n-2}$ containing $I$.

\noindent
Now consider the pull-back ${\bf G}=\{G_i\}_{i=0}^k$ of $F\equiv G_k$ along
the $k$-orbit of $I_0^l$.

{\QP{\bf Lemma 4.2} Under the above circumstances
\item{\rm 1.} $\{G_i\}_{i=1}^k$ is a monotone chain of intervals;
\item{\rm 2.} $G_0\subset T^{l-1}$.\medskip}

{\QP{\bf Lemma 4.3.} The intersection multiplicities of the above chains
${\bf G}$ and ${\bf H}^{n+1}$ are uniformly bounded:
$${\rm mult}\,{\bf G}\leq 8\qquad 
{\rm and}\qquad{\rm mult}\,{\bf H}^{n+1}\leq 8.
$$
\medskip} 

{\bf Proof.} If $t$ intervals of the chain $\{G_i\}_{i=1}^k$ have a common
point, then there is an interval $G_i$ among them containing at least
$(t-1)/2$ intervals $N_s$ of the $(k-1)
$-orbit of $I_1^l.$ Since
$f^{k-i}|G_i$ is monotone, $f^{k-i}N_s$ belongs to the $(u(l-1)-1)$-orbit of
$I_1^l$. But $G_k$ contains at most three intervals of this orbit. Hence
$t\leq 7$.

The argument for ${\bf H}$ is similar, and we omit it.\QED

Now we have enough topological information for getting a priori bounds.

{\QP{\bf Lemma 4.4.}  \qquad$\sup_n \lambda_n\lambda_{n+1}<1.$\medskip}

{\bf Proof.} Choose the smallest interval $I$ among $[0, x_{u(n)}]$ and
$I_k^l\in M^n$ with $k>0$. It is easy to analyse the cases $I=[0, x_{u(n)}]$
or $I=I_k^n$ for $k=1,2.$ So, we restrict ourselves to other cases,
 and then the
interval $F$ is well-defined. Moreover, the Poincar\'e length $[I:F]$ does not
exceed $\log 4$.

It follows from Lemmas 4.2.1, 4.3 and the Schwarz lemma that the Poincar\'e
length $[I_1^l:G_1]$ is uniformly bounded (by a constant depending on 
$\parallel f\parallel$).  Since $f$ is quadratic
 (and hence quasi-symmetric)
near the critical point, the ratio
$${|G_0|\over |G_0\ssm T^l|}$$
can be estimated through $[I_1^l : G_1]$, and hence the ratio 
$|T^l|/|G_0|$ is bounded away from~1.

By Lemma 4.2.2, $G_0\subset T^{l-1}$, so $\lambda_l\leq |T^l|/|G_0|$.
It remains to mention that $\lambda_l$ is equal to either $\lambda_n$
or $\lambda_{n+1}$.\QED

{\QP{\bf Lemma 4.5.} ${1\over{1-\lambda_{n+1}^2}}\leq
      \left({1+\lambda_n\lambda_{n-1}\over 1-\lambda_n\lambda_{n-1}}\right)^2
    (1+O(|{\bf H}^{n+1}|)).$ \medskip} 

{\bf Proof.} Applying the Schwarz lemma to the monotone map
$$f^{u(n)-1} : (H^{n+1},I^{n+1}_1)\rightarrow( T^{n-2} , I_0^n )$$
we get
$$[I_1^{n+1}:H^{n+1}]\leq[T^n:T^{n-2}]+O(|{\bf H}^{n+1}|)
=2\log{1+\lambda_n\lambda_{n-1}\over 1-\lambda_n\lambda_{n-1}}+
O(|{\bf H}^{n+1}|).\eqno (4-2)$$
Let $G$ be the component of $f^{-1}H^{n+1}$ containing 0, $\mu=|T^{n+1}|/|G|.$
The calculation for the quadratic map shows that
$$\log{1\over 1-\mu^2}\leq [I_1^{n+1}:H^{n+1}].\eqno (4-3)$$
Furthermore, since $f^{u(n)}$ is not unimodal on $T^n$,  $G\subset T^n.$
Hence $\lambda_{n+1}\leq\mu .$
The last estimate together with (4-2) and (4-3) yield the required. \QED

From Lemmas 4.4 and 4.5 we get immediately an a priori bound of $\lambda_n$:

{\QP{\bf Lemma 4.6.} $$\sup_n \lambda_n <1.$$\medskip}

{\QP{\bf Lemma 4.7.} Let $L^n$ be the gap between $T^n$ and $J^n$. Then
$$\sup_n{|L^n|\over |x_{u(n)}|}<1.$$
\medskip}

{\bf Proof.} Because of Lemma 4.6, it is enough to show that the gap $L$
is not too small as compared with $J^n$. Let $N$
be a monotonicity  
interval of $f^{u(n-2)}$ adjacent to $J^n$ on its outer side.
Consider the map $f^{u(n-2)}|L\cup J^n\cup N$ and apply to it the Schwarz lemma
taking into account Lemmas 4.1 and 4.6.\QED

Now we can prove that the Lebesgue measure of $M^n$ and ${\bf H}^n$ go down
exponentially fast (compare [G2], [BL3], [MMSS]).
Let $[[\alpha\,,\,\beta]]$ denote the smallest closed interval
containing both $\alpha$ and $\beta$ (similarly, $((\alpha\,,\,\beta))$ 
will denote the smallest open interval containing $\alpha$ and $\beta$).

{\QP{\bf Lemma 4.8.} There exist constants $C>0$ and $q<1$ such that
  $$|{\bf H}^n|\leq Cq^n\qquad {\rm and} \qquad |M^n|\leq Cq^n.$$
Hence, the Lebesgue measure of $\omega(0)$ is equal to zero.\medskip}

{\bf Remark.} The last statement is a corollary of more general
results [BL2],~[M].

{\bf Proof.} By Lemma 4.7, density of $M^{n+1}$ in $I_0^n$ is bounded away
from 1.
Consider now an interval $I^n_l\in M^n,\;  l>0.$ It follows from
Lemmas 4.1 and 4.6 that the map
$$f^{u(n-1)-l}: I_l^n\rightarrow [[x_{u(n-1)}, x_{u(n+1)}]]$$
has  bounded distortion. But this map carries $M^{n+1}\cap I^n_l$ into
$I_0^{n+1}\cup J^n.$ By Lemma 4.7, density of the latter set in 
$[[x_{u(n-1)}, x_{u(n+1)}]]$ is
bounded away from 1. Hence density of $M^{n+1}$ in $I^n_l$ is bounded away
from 1 as well. So, there is a $q<1$ such that
$$\lambda(M^{n+1})\leq q\lambda(\cup_{l=0}^{u(n-1)-1}I^n_l)+
\lambda(\cup_{l=u(n-1)}^{u(n)-1}I_l^{n+1}).$$
Applying this estimate twice we get
$$\lambda(M^{n+2})\leq q\lambda(M^n),$$
and we are done with $M^n$.

Now consider a pair $H^{n+1}\supset H^{n+2}$ and apply $f^{u(n)-1}.$ Then
$H^{n+1}$ is mapped onto $T^{n-2},$ while $H^{n+2}$
is mapped into $T^{n-1}$ (since $f^{u(n-1)}$ is monotone on its
image). By Lemma 4.6 and the Schwarz lemma, the density of $f^m H^{n+2}$ in
$f^m H^{n+1}$ is bouded away from 1 for $m=0,...,u(n)-1.$ Furthermore,
$$f^{u(n)+m} H^{n+2}\subset I_m^{n-1}\; ,\qquad m=1,...,u(n-1).$$ 
Cosequently, for some $q_1<1$ we have
$$|{\bf H}^{n+2}|\leq q_1 |{\bf H}^{n+1}|+|M^{n-1}|+|M^{n-2}|,$$
and the required follows.\QED\smallskip

{\QP{\bf Lemma 4.9.} (i). There is a $q<1$ such that
 $\lambda_{n+1}^2=O(\lambda_n\lambda_{n-1}+q^n).$
{\item (ii)}. $\lambda_{n+1}^2=O\left({|J^n|\over |T^{n-1}|}\right).$ 
\medskip}

{\bf Proof.} The point (i) follows from Lemmas 4.5 and 4.8.
To prove (ii), consider $f^{u(n-1)}: I^{n+1}\rightarrow J^n$
and apply the Schwarz lemma.\QED

{\bf Remark 4.10.} All constants in the above estimates depend only on 
$\parallel f\parallel.$ Moreover, they are uniform over the maps with negative
Schwarzian derivative (since the Schwarz lemma and the Koebe Principle are
uniform over this class). Finally, all estimates are asymptotically uniform
over the whole class ${\cal U}$ ("beau estimates", see Sullivan [S]).
 For example,
Lemma 4.6 can be improved in such a way:
$$\limsup_{n\to\infty}\lambda_n\leq C<1$$
for an absolute constant $C$.\bigskip 

\centerline{\bf \S5 Scaling, characteristic exponent and Hausdorff dimension.}
\bigskip

In this section we will prove Theorem 1.3 assuming that there is a good enough
a priori bound of $\lambda_n$. It follows that the Theorem holds for an open set
of Fibonacci maps invariant under quasi-symmetrical conjugacy.

Let $q<1$ be the constant from Lemma 4.8 , 
$\sigma_n=\max_{n-1\leq i\leq n+1}(\lambda_i,\lambda_{i+1})$.

{\QP{\bf Lemma 5.1.} For any $x\in I_1^{n+1}$
$${d_n\over d_{n+1}^2}(1+O(\sigma_n+q^n))^{-1}\leq |(f^{u(n)-1})'(x)|\leq
{d_n\over d_{n+1}^2}(1+O(\sigma_n+q^n)).$$\medskip}

{\bf Proof.} Let us apply the Koebe Principle  to the map 
$$f^{u(n)-1}: (H^{n+1}, I^{n+1})\rightarrow (T^{n-2}, T^n)$$ taking into
account Lemma 4.8: 
 $${|(f^{u(n)-1})'(x)|\over |(f^{u(n)-1})'(y)|}=1+O(\lambda_n\lambda_{n-1}+q^n),
\qquad x, y\in I_1^{n+1}.$$
Besides,
$${d_n\over d_{n+1}^2}\leq{|I^n|\over |I_1^{n+1}|}\leq 
(1+\lambda_{n+1}\lambda_{n+2}){d_n\over d_{n+1}^2},$$
and the Lemma follows.\QED
 
{\QP{\bf Lemma 5.2.} There is a $\rho=\rho(\parallel f\parallel)$ and 
$L=L(\parallel f\parallel)\in {\bf N}$
such that if $\lambda_l <\rho$ for some $l\geq L$ then $\lambda_n$ exponentially
decrease. For maps with non-positive Schwarzian derivative one can
 choose $L=1$ and uniform $\rho$. \medskip}

{\bf Proof.} Let $n$ be so large that $f(x)$ is a quadratic map in the 
neighborhood $T_{n-1}$. Then by the chain rule,
$$|(f^{u(n)-1})'(x_1)|=|(f^{u(n-1)-1})'(x_1)|\cdot 2 d_{n-1}
|(f^{u(n-2)-1})'(x_{u(n-1)+1})|.\eqno (5-1)$$
By Lemma 5.1,
$${d_{n-1}\over d_n^2}\cdot 2 d_{n-1}{d_{n-2}\over d_{n-1}^2}\leq
{d_n\over d_{n+1}^2}(1+O(\sigma_{n-2}+\sigma_{n-1}+\sigma_n+q^n)).\eqno (5-2)$$
It follows from Lemma 4.5 that $\lambda_k$ keep to be small for $k=n-2,...,n+1$,
once $\lambda_{n-3}$ becomes small for big enough $n$. Hence, by (5-2)
$$\lambda_{n+1}^2\leq\gamma\lambda_n\lambda_{n-1}\eqno (5-3)$$
for some $\gamma<1$.
Setting $\Lambda_n=\max(\lambda_n , \lambda_{n-1}),$ we get from (5-3) that
$$\Lambda_{n+1}\leq\sqrt{\gamma}\Lambda_n.\eqno (5-4)$$
So, once  $\lambda_n$ become small, they start exponentially decrease.
It follows that they exponentially decrease forever.

The final remark: since the constants in the Schwarz Lemma and the Koebe 
Principle depend only on $\parallel f\parallel$
, the constants $\rho$ and $L$ depend only on this data as well .  Moreover,
all estimates are  uniform in the case of negative Schwarzian derivative.\QED

Recall that a one dimensional homeomorphism  $h$ is called {\sl quasi-symmetric}
if any two adjacent commensurable intervals $I$ and $J$ are mapped into
commensurable ones:
$${|I|\over|J|}\leq K\Rightarrow {|fI|\over |fJ|}\leq \gamma(K).$$

Denote by ${\cal F}^0$ the set of Fibonacci maps for which $\inf \lambda_n=0$.
{\QP{\bf Lemma 5.3.} 
\item{\rm 1.} The set ${\cal F}^0$ is invariant under quasi-symmetrical
conjugacy.\smallskip
\item{\rm 2.} The set ${\cal F}^0$ is $C^0-$open in the $C^2$-balls $B(r)$
of the space ${\cal F}$.\medskip}

{\bf Proof.} The first point is clear from the definitions .
 The second one follows from the fact that the constants in the previous 
lemma are uniform over $B(r)$.\QED

Let us write $\alpha_n\sim\beta_n$ if $|\log{(\alpha_n/\beta_n)}|$ 
exponentially decrease, and  $\alpha_n\asymp\beta_n$ if it is bounded.

The next lemma gives the asymptotical formula of Theorem 1.3.1 for the subclass~
${\cal F}^0$ (compare Tangerman and Veerman [TV]).

{\QP{\bf Lemma 5.4.} For any $f\in{\cal F}^0$ the following asymptotical formulas
hold:
\item{\rm 1.} $\lambda_{n+1}\;\sim\;\lambda_n/{\root 3\of 2}.$
\item{\rm 2.} $\lambda_n\;\sim \;a 2^{-n/3}.$
\item{\rm 3.} $d_n\;\sim\;(1/2)^{n^2/6+\beta n+\gamma}$

for some constants $a>0 , \;\beta $ and $\gamma$.
Moreover
$$|\log (a/\lambda_0)|
\leq R(\parallel f\parallel)\;,$$
and the constant $R$  is uniform over 
maps with negative Schwarzian derivative.\medskip}


{\bf Proof.} Since $\lambda_n$ exponentially decrease, Lemma 5.1 yields for
${x\in I_1^{n+1}}$
$$|(f^{u(n)-1})'(x)|\;\sim\; d_n/d_{n+1}^2.\eqno (5-5)$$
Substituting this into the recurrent equation (5-1), we get
$$\lambda_{n+1}^2\;\sim\;{1\over 2}\lambda_n\lambda_{n-1}.\eqno (5-6)$$
Setting $s_n=\log(\lambda_n/\lambda_{n-1})$, we have from the last formula
$$s_{n+1}=-{1\over 2}s_n-{1\over 2}\log 2 +O(q^n) $$
with $q<1$. It yields
$$s_n=-{1\over 3}\log 2+O(\rho^n)\eqno (5-7)$$
with $\rho=\max(1/2, \; q)$ which proves the first point of the lemma.

Setting now $c={1\over 3}\log 2, \; \nu_n=\log \lambda_n+cn$ we get from (5-7)
$$\nu_{n+1}=\nu_n+O(\rho^n).$$
So, there is a limit 
$$\lim\nu_n\equiv\log a=\nu_0+O(1),$$
with exponential convergence
and the constant depending only on 
$\parallel f\parallel$ and uniform over maps with negative Schwarzian. 
Equivalently
$$a\equiv\lim\lambda_n e^{nc}\asymp\lambda_0. $$
It proves the second point together with the last remark.
 The reader can easily derive the third point from the second one.\QED

Let us estimate now the ratio of any two intervals 
$M^n_{s_1...s_n}\subset M^{n-1}_{s_1...s_{n-1}}$. 
  The previous lemma gives the asymptotics for the ratio
$\lambda_n\equiv |M^n_{0...0}|/|M^{n-1}_{0...0}|$ . Besides,
$M^n_{s_1...10}=M^{n-1}_{s_1...1}$. Other cases are covered by the following
 lemma.

{\QP{\bf Lemma 5.5.} For $f\in {\cal F}^0$ the following scaling laws hold:
$${|M^n_{0...01}|\over |M^{n-1}_{0...0}|}\;\equiv\;{|J^n|\over |I^{n-1}|}\;\sim\;
{a^2\over 2^{2(n+1)/3}}.$$
If $[s_1...s_{n-1}]\neq [0...0]$ then
$${|M^n_{s_1...s_{n-1}1|}\over |M^{n-1}_{s_1...s_{n-1}}|}\;\sim\;
            {a^2\over 2^{2(n-1)/3}}.$$
and
$${|M^n_{s_1...s_{n-1}0|}\over |M^{n-1}_{s_1...s_{n-1}}|}\;\sim\;
            {a^2\over 2^{2(n-2)/3}}$$
where $a$ is the constant from Lemma 5.4. All asymptotics are uniformly
exponential.\medskip}

{\bf Proof.} Let us consider a chain of two maps

$$\matrix
{(I^{n-1}, J^n)&\rightarrow &(I_1^{n-1}, J_1^n)&\rightarrow &(I^{n-2}, I^n).\cr
 &f&&f^{u(n-2)-1}&\cr}$$
Note that by Lemma 5.4 $|I^n|\sim |x_{u(n)}|$.
Setting $r_n=|J^n|/|I^{n-1}|$ we get
$${|fJ^n|\over |fI^{n-1}|}\sim 1-(1-r_n)^2\;\sim\;2 r_n.$$
On the other hand, $f^{u(n-2)-1}$ has an exponentially small distortion on
$I_1^{n-1}$ , and hence
$$2r_n\;\sim\;{|I^n|\over |I^{n-2}|}\sim \lambda_n\lambda_{n-1}\;\sim\;
{a^2\over2^{(2n-1)/3}},$$
and the first asymptotical formula is proved.

In order to get the others, consider the map
$$f^k: \;M^{n-1}_{s_1...s_{n-1}}\rightarrow I^{n-2}$$
for an appropriate $k$. It carries $M^n_{s_1...s_{n-1}0}$ into $J^{n-1}$
and $M^n_{s_1...s_{n-1}1}$ into $I^n$ with exponentially small distortion.
It yields the required.\QED

Now we can prove the next piece of Theorem 1.3 for $f\in {\cal F}^0$

{\QP{\bf Lemma 5.6.} For $f\in{\cal F}^0 $ the critical orbit closure 
$\overline{\cal O}$
has Hausdorff dimension~0.\medskip}
 {\bf Proof.} Let us consider covering of $\overline{\cal O}$ by the intervals
$M^n_{s_1...s_n}$. By the above two lemmas, the lengths of these intervals
decrease uniformly superexponential\break
($O(q^n) $ for any $q\in (0,1)$), while their number increases exponentially
 ($\leq 2^n$).  Let $\gamma=-\log 2/\log q ,\qquad l_\gamma$ be the Hausdorff
measure on $\overline{\cal O}$ of exponent $\gamma$. Then
$$l_\gamma(\overline{\cal O})\leq C 2^n q^{n\gamma}\leq C.$$
Hence, $\dim \overline{\cal O}\leq\gamma$, and $\gamma$ is arbitrary small
positive number. \QED

Now we are going to show that the geometry of the set $\overline{\cal O}$
is completely determined by only one parameter $a$ from Lemma 5.4. Let $f$ and
$g$ be two Fibonacci maps,
$$\phi : \overline{\cal O}(f)\rightarrow\overline{\cal O}(g)$$
be the natural topological conjugacy. Let us say that $\phi$ is smooth if
for any $x\in\overline{\cal O}$ there exist
$$\lim {|\phi(x)-\phi(y)|\over |x-y|}\neq 0$$
as $y\to x$ along $\overline{\cal O}(f)$, 
and this limit depends continuously on $x$.

{\QP{\bf Lemma 5.7.} If two Fibonacci maps $f$ and $g$ in ${\cal F}^0$  have the 
same parameter $a$ then the conjugacy $\phi$ is smooth on $\overline{\cal O}(f)$.
\medskip}  

{\bf Proof.} Indeed, it follows from Lemmas 5.4 and 5.5 that for any Fibonacci
sequence ${\overline s}=s_0 s_1...$ there is a limit 
$$\lim_{n\to\infty}{|M^n_{s_1...s_n}(f)|\over |M^n_{s_1...s_n}(g)|} $$
depending continuously on ${\overline s}$.\QED

{\QP{\bf Lemma 5.8.} Let $f\in {\cal F}^0 ,\qquad n=[s_1...s_k]$ be the Fibonacci
expantion of $n$. Then
$$|(f^n)'(x_1)|\;\sim\;2^{{2\over 3}\sum{ms_m}+\gamma\sum{s_m}+\delta}$$
for some constants $\gamma$ and $\delta$.\medskip}

{\bf Proof.} Let $m_i$ be the places where$s_{m_i}=1$. Decompose $n$-orbit
of $x_1$ into the parts of length $u(m_i)$. By (5-5) it gives the factorization
of the derivative into factors of order $\sim 2\lambda^{-2}_{m_i+1}$. Now
Lemma 5.4 implies the required asymptotics.\QED 

Clearly, it follows from the last lemma that the growth of the $n$-fold derivative at $x_1$ is subexponential. The maximal growth of order 
 $\exp\kappa(\log n)^2$  (which is faster than any power $n^\gamma$)
is attained at noments $u(m)-1$. However,
at the next moments $n=u(m)$ the derivative drops to $n^\gamma$ with 
$\gamma=2\log 2/ 3\log ({\sqrt{5}+1\over 2})<1.$ These oscillations are balanced
in a  ``convergent way" .

{\QP{\bf Lemma 5.9.} The series
$$\sum_{n=1}^\infty{1\over |(f^n)'(x_1)|^\alpha}$$
 is convergent for any $\alpha>0$.\medskip}

{\bf Proof.} By the last lemma , this series has a majorant of the following
form:
$$\sum_{s_m\in\{0,1\}} 2^{-\sum_{m=1}^k(am+b)s_m}=
\prod_{m=1}^\infty (1+{1\over 2^{am+b}})<\infty.$$\QED 

This Lemma and the Nowicki-van Strien Theorem [NvS] imply the
existence of an
absolutely continuous invariant measure for $f\in {\cal F}^0$. So, Theorem 1.3
is proved for the subclass ${\cal F}^0$.
\bigskip
\eject
\centerline{\bf \S6. Real renormalizations.}\medskip

Now we need another class of maps on which we can define a renormalization in
such a way that the Fibonacci maps can be exactly characterized as
infinitely renormalizable. Let
$$	J=[a,b]\;,\qquad T=[\alpha , \beta]\;,\; {\rm where}\;
-1<a<b<\alpha<\beta<1\;,\qquad Dom(f)=J\cup T\;,$$
and let $\;f: Dom(f) \rightarrow [-1 , 1]\;$
be a $C^2$-smooth map such that (see Figure 2) :

\item{(i)} $ f|J$ is a diffeomorphism from $J$ onto $[-1,1]$, which may
be either orientation preserving or orientation reversing.

\item{(ii)} $f|T$ is a unimodal map from $T$ into $[-1,1]$ with
non-degenerate 
minimum point, and with $f(\partial T)=1$.

\noindent
 Let us denote the space of all such maps by ${\cal A}$. 
Since we don't specify whether $f|J$ preserves or reverses orientation,
 ${\cal A}$ can be decomposed into the union of two connected
components ${\cal A}^+$ and ${\cal A}^-$, where ``+" corresponds to the case of
orientation preserving $f|J$. \smallskip

\centerline{\psfig{figure=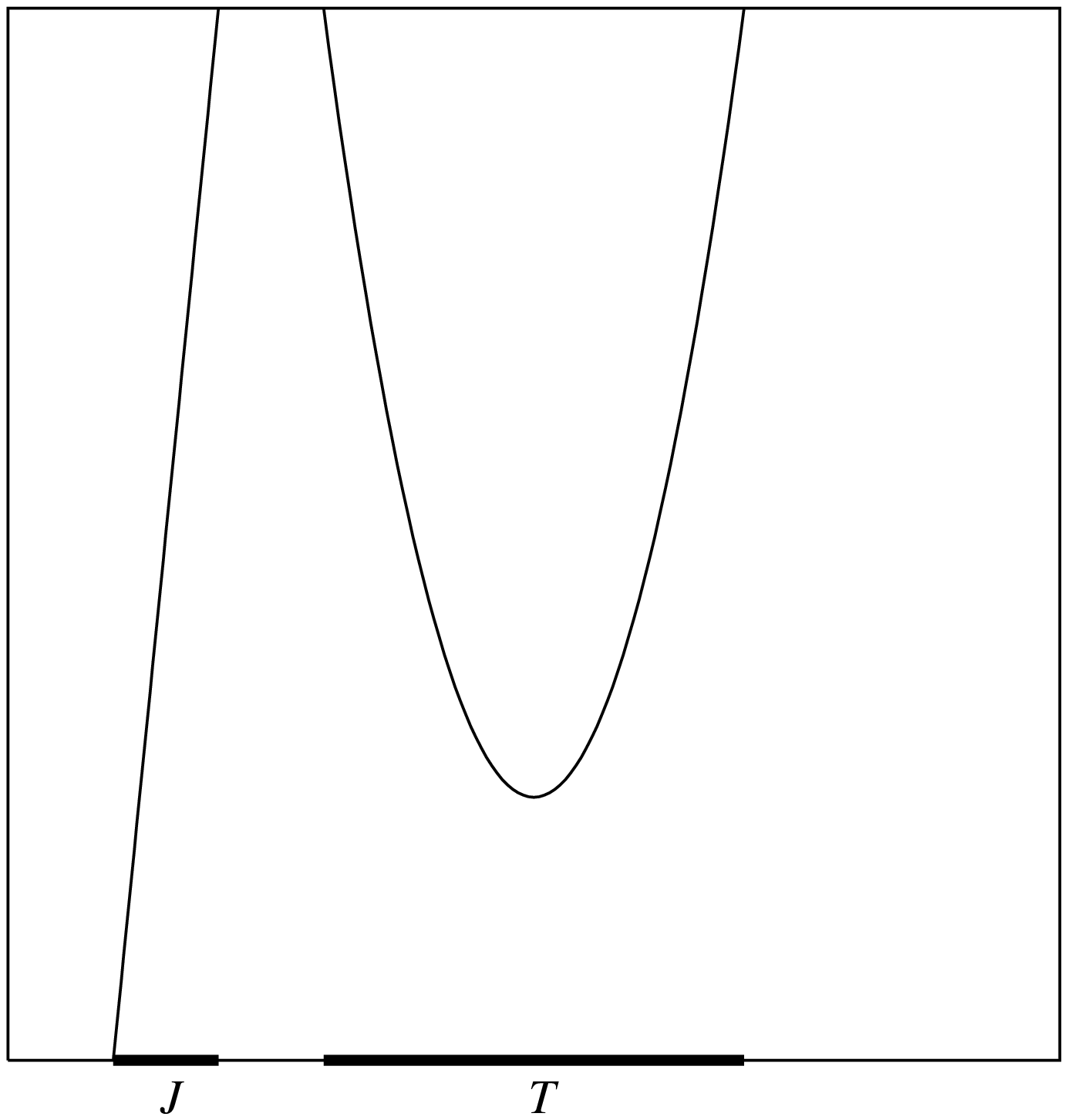,height=2in}}
\centerline{Figure 2, Graph of a function in ${\cal A}^+_0$.}\smallskip

Now suppose that some map
$f\in{\cal A}$, with critical point $x_0\in T$, satisfies the conditions
that the critical value $x_1=f(x_0)$ lies in $J$, and the its image
$x_2=f(x_1)$ lies back in $T$. Then we will be interested in two segments
of the first return map from $T$ to itself, as follows. There is an interval
$T_1$ around the critcal
point which is mapped unimodally by $f^2$ into $T$, with both endpoints
of $T_1$ mapping to one endpoint of $T$. 
Further, there is a disjoint interval $J_1\subset T$
which maps diffeomorphically onto $T$ under the map $f$ itself. Here we choose
$J_1$ to the left of $T_1$ if $f|J$ preserves orientation, or to the right
of $T_1$ if $f|J$ reverses orientation
(so that $J_1$ lies on the same side of 0 as $x_2$).
 The resulting map $Vf : J_1\cup T_1\to T$,
affinely conjugated ({\sl rescaled}) 
so that $T$ is replaced by the original interval
$[-1,1]$, is the required renormalization $Rf$
(there is choice of two rescalings ; select that one which makes
 the critical point
to be minimum point ). This renormalization interchanges
the two spaces ${\cal A}^+$ and ${\cal A}^-$. If $f$ is $n$-fold
renormalizable then $R^nf$ comes as rescaling of a map $V^nf\equiv f_n$,
the restriction of appropriate iterates of $f$ to the union of two
appropriate intervals, $T_n$ and $J_n$.

Let $T_+$ and $T_-$ be the semi-intervals on which 0 divides $T$.
The {\it kneading sequence} of $f\in {\cal A}$ is the sequence of symbols
$U_n\in \{T_+,T_-,J\}$ such that $x_n\equiv f^n0\in U_n$. Two maps
$f\in {\cal A}^+$ (or ${\cal A}^-$) without limit cycles are topologically
conjugate if and only if they have the same kneading sequence (compare [MT]).

In terms of kneading sequences the above renormalization can be described
in the following way. The renormalizable kneading sequences start with $JT_s,\; 
s\in\{+,-\} $. To write its renormalization do the
following operations moving along the sequence:
{\item (i).} When you see $J$, cross it;
{\item (ii).} When you see $T_sJ,\; s\in \{+,-\}$, change $T_s$ for $T_{ks}$
provided $f\in {\cal A}^k,\; k\in\{+,-\}$.
{\item (iii).} When you see $T_sT_r$, change the first $T_s$ for $J$.

Let us say that a map $f\in {\cal A}^+$ is a Fibonacci map if it has the
following kneading sequence:
$$fib^+=J|T_-|T_+|JT_+|JT_-T_-|JT_-T_+JT_-|...$$
(In order to write the block from $u(n)+1$ to $u(n+1)$ repeat the beginning
of the sequence till the moment $u(n-1)$, and then change the last symbol
$T_s$ for the ``opposite" one, $T_{-s}$). Denote this class of maps by
${\cal F}^+$. Similarly, the kneading sequence of a map $f\in {\cal F}^-$
is produced by the same rule but with different initial:
$$fib^-=J|T_+|T_+|JT_-|JT_+T_-|JT_+T_+JT_+|...$$
A class ${\cal F}$ of Fibonacci maps is defined as 
${\cal F}^+\cup{\cal F}^-$. One can also describe this class 
by the following properties: $x_1\in J,$ and $f^{u(n-1)}$ is well-defined
and monotonous on the interval $[[0,x_{u(n)}]]$, and 
$$f^{u(n-1)}((0,x_{u(n)}))\equiv ((x_{u(n-1)},x_{u(n+1)}))\ni 0. \eqno (6-0)$$

   If we want to emphasize that $f\in {\cal A}$ 
then we say that $f$ has type
(2,1). In the unimodal case  we say that $f$ is of type (2)
(see the next section for more general discussion). As in the unimodal
case, we will use the notations $T^n$ and $J^n$ for the intervals
$[[x_{u(n)},x_{u(n)}']]$ and $[[x_{u(n-1)}, x_{u(n-1)+u(n+1)}]]$
correspondingly (don't confuse with $T_n$ and $J_n$ introduced above).
  
{\QP{\bf Lemma 6.1.} A map $f\in{\cal A}$ is infinitely renormalizable if and
only if it is  a Fibonacci map : $f\in {\cal F}$.
In this case the following inclusions hold:
$$T^{n+2}\subset T_n\subset T^{n+1} \eqno (6-1) $$
$$J^{n+2}\subset J_n. \eqno (6-2)$$\medskip}

{\bf Proof.} Let $f\in {\cal A}$ be infinitely renormalizable.
 Then one can check by induction that\smallskip
   $$f_{n}|T_{n}=f^{u(n+1)}\qquad {\rm and} \qquad f_n|J_n=f^{u(n)}.\eqno (6-3)$$
Since $f_{n-1}$ is renormalizable,   
$$x_{u(n)}=f_{n-1}(0)\in J_{n-1}\qquad {\rm and} \qquad
x_{u(n+1)}=f_n(0)\in T_{n-1}.$$
Hence, $x_{u(n+1)}$ lies closer to 0 than $x_{u(n)},\; n=1,2,...$

Let us study now the combinatorics of several first iterates of 0. Since $f$
is renormalizable, 
$$T^2\equiv[[x_2,x_2']]\subset T\subset [x_1,x_1']\equiv T^1. \eqno (6-4)$$
Furthemore, $x_3=f_2(0)\in T_1$; hence $x_4=fx_3\in J.$ So,
$$ J^2\equiv[x_1,x_4]\subset J.\eqno (6-5)$$

 Consider now the following map 
$\sigma: {\bf N}\rightarrow {\bf N}$ of the set of natural numbers: if 
$m=\sum u(l_i)$ is the Fibonacci expantion of $m$ then
$\sigma(m)=\sum u(l_i+1)$ ($\sigma$ is induced by the shift on the 
space of Fibonacci expantions). Then we have the following rule:
    $$(f_n)^m(0)=x_{\sigma^n(m)}. \eqno (6-6)$$
So, if we have a combinatorial property of several points $x_m$ then repalcing
$f$ by $f_n$ we immediately get the same property of points $x_{\sigma^n m}$
(provided $f$ is infinitely renormalizable). In particular we can replace
ponts $x_1, x_2, x_4$ in (6-4) and (6-5) by\break
$x_{u(n+1)}, x_{u(n+2)}, x_{u(n+1)+u(n+3)}$.  
Then we obtain the required properties (6-1) and (6-2).

Let us show now that $x_1$ and $x_2$ lie on the same side of 0 for
$f\in{\cal A}^+$, and they lie on the opposite sides of 0 for $f\in{\cal A}^-$.
Indeed, otherwise consider $f|[x_1, x_4]$ and conclude that $x_5$ lies
farther from 0 than $x_2$.

Changing $f$ for $f_1$ we get the same statement for the points $x_2$ and
$x_3$. Since the renormalization interchanges ${\cal A}^+$ and
${\cal A}^-$, we conclude that $((x_1,x_3))\ni 0$. Replacing $f$ by $f_{n-2}$
we get (6-0).

Finally, since $x_2\in T,\; f | [0,x_2]$ is well-defined and monotone. 
Replacing it again by $f_{n-2}$ we conclude that 
$f^{u(n-1)} | [0,x_{u(n)}]$ is well-defined and monotone. So, $f$ is a
Fibonacci map.

Vice versa, let $fib^s_n,\; s\in \{+,-\}$ be the initial parts of
length $u(n)$ of the kneading sequences $fib^s$. Then one can easily
check by induction that the renormalization turns $fib_n^s$ into
$fib_{n-1}^{-s}$. So, it interchanges $fib^s$ and $fib^{-s}$ which certainly
implies that both sequences are infinitely renormalizable. 
\QED 

Now let us briefly discuss topology on the space $\cal A$ (compare \S 4).
 We can restrict ouselves to  the subspace ${\cal A}_0\subset{\cal A}$
consisting of those $f$ for
which $f|T$ is an even function, $f(-x)=f(x)$.   
 Then we can
write $f|T$ uniquely as $$f(x)=A_{x_1}\circ f_T\circ Q\circ A_T$$ where $A_T$ is
the orientation preserving linear map which carries $T$ onto $[-1,1]$, $Q$
is the squaring map
$\xi\mapsto\xi^2$, $f_T$ is some orientation preserving diffeomorphism
of $[0,1]$, and $A_{x_1}$ is the orientation preserving affine map which carries
$[0,1]$ onto $[x_1,1]$, where $x_1=f(0)$ is the critical value. Similarly, we
can write $f|J$ as $f_J\circ A_J$ where $A_J$ is the orientation preserving
affine map from $J$ onto $[-1,1]$, and where $f_J$ is a diffeomorphism
of $[-1,1]$.

Now we suppose that both $f_J$ and
$f_T$ are $C^2$-smooth. 
The  $C^k$-topology on ${\cal A}_0,\; k\leq 2,$
 comes from the $C^k$-topology on the space
of diffeomorphisms $f_T$ and $f_J$, together with the Euclidian topology on
the finite
dimensional space of parameters $a,b,\alpha,\beta, x_1$. 
Let $\| f\|$ denote the maximum of the $C^2-$norms for
$f_J , f_J^{-1}$ and $f_T , f_T^{-1}$ which is a continuous functional on our
 space. 

We can  assume without loss of generality that the original map
$f$ is quadratic near 0 (though this property is not preserved under
renormalization). Let us remark also that clearly all estimates of
\S\S 4,5 hold not only for unimodal maps but in the class ${\cal A}$
as well.

{\QP{\bf Lemma 6.2.}  The norms $\|R^nf\|$  
are uniformly bounded.\medskip}

{\bf Proof.} By (6-3), $f_n|T_n=f^{u(n+1)}$ which can be decomposed as a
quadratic map and the diffeomorphism
$$f^{u(n+1)-1} : H^{n+2}\rightarrow T^{n-1} \eqno (6-6)$$
(see Lemma 4.1). On the other hand,
$$f^{u(n+1)-1}(fT_n)=f_nT_n\subset T_{n-1}\subset T^n  \eqno (6-7)$$
(the last inclusion is by (6-1)). It follows from (6-6), (6-7) and
a priori bounds proven in \S 4 that $f^{u(n+1)-1}| fT_n$ has bounded
distortion. By rescaling we get
$$\log \left|{(R^nf)_T'(x)\over (R^nf)_T'(y)}\right|=O(|x-y|)$$
for any $x, y\in [0,1].$ This implies
$$\left|{(R^nf)_T''\over (R^nf)_T'}\right|=O(1).$$

 Because of bounded
distortion, the derivative $(R^nf)_T'$ is uniformly bounded from 
below and above, and  the boundedness property for
the second derivaty $(R^nf)_T''$ follows.  The same argument applies
to $(R^nf)_J$ and to the inverse maps.\QED

{\bf Corollary.} If $\inf \lambda_n>0$ then there is a $C^1$-convergent sequence
of renormalizations $R^{n_i}\to g\in {\cal A}$.

{\bf Proof.} It follows from the assumtion and inclusions (6-1) that the ratio 
$|T_n|:|T_{n-1}|$ is bounded away from 0. Moreover, Lemma 4.9 and (6-2)
imply the same for the ratio $|J_n| : |T_{n-1}|$.
Now one can play the ``distortion game" in manner of \S 4 to check that
three complementary gaps (that is, components of $T_{n-1}\ssm(T_n\cup J_n)$)
are also commensurable with $T_{n-1}$. After rescaling we conclude that
the domains $Dom(R^nf)$ don't degenerate, so we can select a convergent
sequence $Dom(R^{n_i}f)$.  Then by the
last lemma,  families of diffeomorpfisms $\{(R^{n_i}f)_T\}$ and $\{(R^{n_i}f)_J\}$
are $C^1$-precompact, and we can extract from them convergent subsequences
as well. \QED

For an interval $I\subset {\bf R}$
denote by $P(I)$ the plane slitted along two rays: 
$$P(I)={\bf C}\ssm({\bf R}\ssm I). $$

Let us introduce now a subspace ${\cal E}\subset {\cal A}$ consisting of maps 
$f: T\cup J\rightarrow [-1,1]$ with the following property:
The map $f_T^{-1}: [0,1]\rightarrow [0,1]$ can be analytically continued
to a map $P[0,1]\rightarrow P[0,1]$,
  and $f^{-1}: [-1,1]\rightarrow J$ can be analytically continued
to a map $P[-1,1]\rightarrow P(J)$.
   
{\QP{\bf Lemma 6.3.} Let $ R^{n_i}f\to g$ in $C^1$-topology. Then the limiting
function $g$ belongs to the class ${\cal E}$.\medskip}

{\bf Proof.} The map $(R^nf)_T^{-1}$  can be written
as long compositions of type $h_1\circ q_1\circ...\circ h_k\circ q_k$
where $h_i$ are diffeomorphisms between apropriate intervals
with a small total distortion while
$q_i$ are square root maps (we reserve this term for affine conjugates to
the standard square root). Such a map can be rewritten as $H_n\circ Q_n$
where the distortion of $H_n$ does not exceed the total distortion of 
$h_i,\; i=1,...,n$,
and $Q_n$ is the composition of $Q_i$ renormalized by appropriate M\"{o}bius maps
(see [S] , [Sw2]).  The maps $Q_n$ analytically map $P[-1,1]$ into itself,
and hence form a normal family. So, we can select a convergent sequence
$Q_n\to Q$ with $Q$ to be a self-map of $P[-1,1]$. On the other hand,
$H_n\to$id in $C^1$-topology. So, $(g_T)^{-1}=Q$. In the same way we can 
treat $g_J$. \QED

{\bf Correspondence between Fibonacci maps of classes ${\cal U}$ and 
${\cal A}^-$.}
We are going to describe an easy surgery interchanging these classes 
without touching the critical orbit. It will follow that any result about
the critical orbit established in one of the classes immediately yields
the same statement in the other class.

Let $f\in {\cal U}$ be a unimodal Fibonacci map. Let us restrict it onto the
union of two disjoint intervals 
$$I^2\cup J^2\equiv [x_5, x_2]\cup [x_1, x_4].\eqno (6-8)$$
Then let us embed these intervals into disjoint intervals $T$ and $J$
correspondingly, and  continue $f$ to a map of class ${\cal A}^-$ defined on
$T\cup J$.  

Vice versa, given a Fibonacci map $g\in{\cal A}^-$,  we can also restrict it
onto the union (6-8), and then continue to a unimodal map of class 
${\cal U}$. This is possible since 
$g(x_5)\equiv x_6<x_5\equiv g(x_4)$.

Since  orb$(0)\subset I^2\cup J^2$, the above surgeries keep the critical
orbit untouched.

\bigskip 
\centerline{\bf \S 7. Polynomial-like maps.}\medskip

Now we are going to show that all polynomial-like maps $f\in {\cal A}^-$
(or ${\cal A}^+$)
are quasi-symmetric\-ally conjugate. It is convenient to introduce more
general terminology. 

Consider $k+1$ topological disks $U_i$ and $V$
bounded by piecewise smooth curves, and such that 
${\rm cl} U_i$ are disjoint and contained in $V$. 
Let us say that
$$f : \cup U_i\rightarrow V$$
is {\it a polynomial-like map of type} $( n_1,...,n_k)$   
if $f|U_i$ is a branched covering of degree $n_i$; $d=\sum n_i$
is called the degree of $f$. 
 Note that polynomial-like maps of type $(d)$ are exactly polynomial-
like maps in the sense of Douady and Hubbard [DH]. 

{\QP{\bf Lemma 7.1.} Any  polynomial-like map 
$f:U_1\cup U_2\rightarrow V$ of type $(2,1)$ is 
quasi-conformally conjugate to
a cubic polynomial with at least one escaping critical point.\medskip}

{\bf Proof.} Consider an ``eight-like" neighborhood $N$ of $U_1\cup U_2$ 
and smoothly continue $f$ there so that $f$ becomes a double covering on the
annulus around $U_1$ and a diffeomorphism on the annulus around $U_2$, and
both annuli are mapped on the same annulus around $V$, see Figure 3.

Then continue $f$ to a slightly bigger domain so that it turns into a three
sheeted smooth covering of a topological disk over a bigger disk.  Now use
the Douady-Hubbard surgery [DH] in order to quasi-conformally conjugate
this map to a cubic polynomial.\QED

\centerline{\psfig{figure=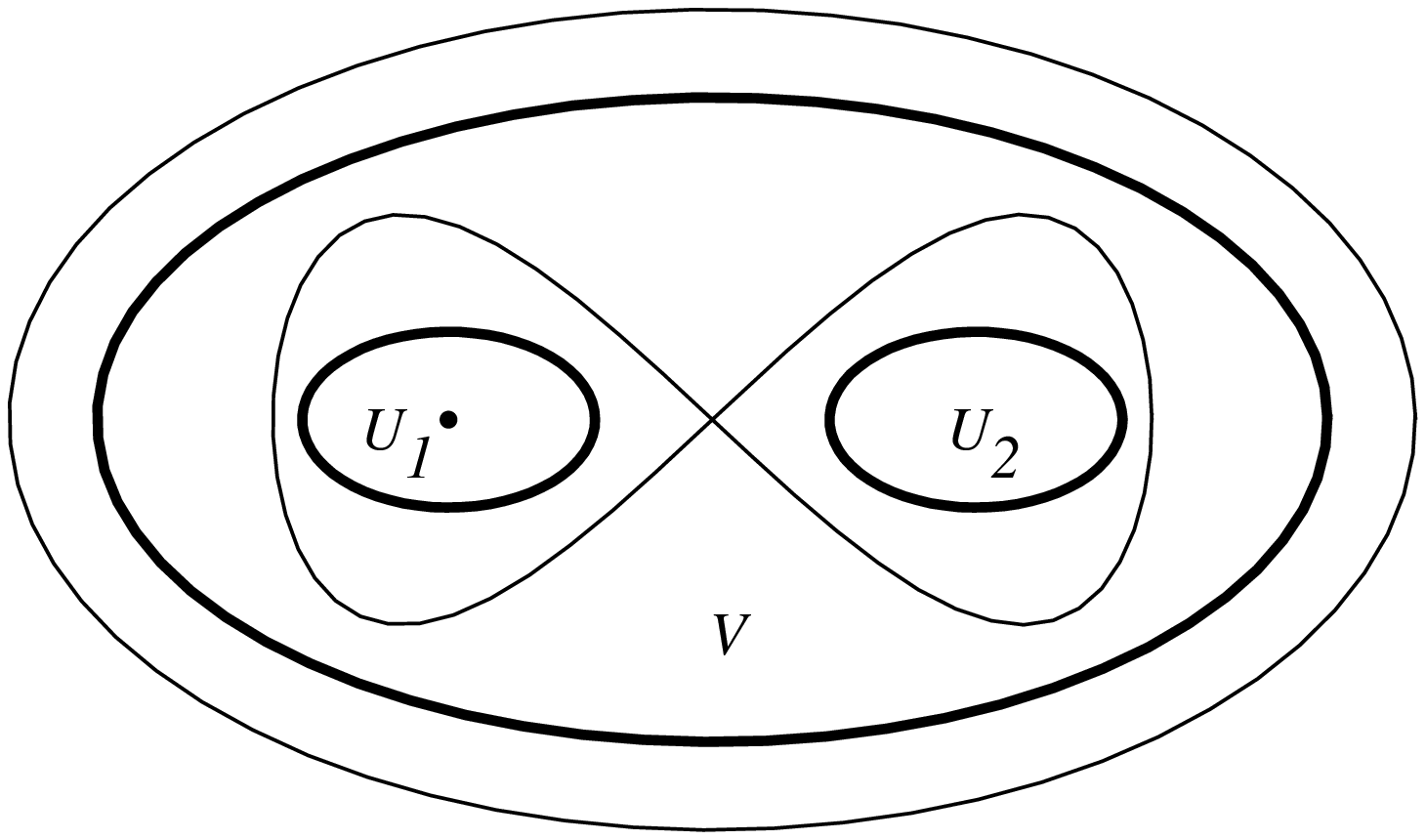,height=2in}}
\centerline{Figure 3.}
\bigskip

{\QP{\bf Lemma 7.2.} Any polynomial-like map $f\in {\cal A}^-$ is 
quasi-symmetrically conjugate to a real cubic polynomial with
one escaping critical point.\medskip}

{\bf Proof.} For $f\in {\cal A}^-$ one can carry out the above construction
in an $\bf R$-symmetrical way. \QED

{\QP{\bf Lemma 7.3.} All Fibonacci real cubic polynomials are quasi-symmetrically
conjugate.\medskip}

{\bf Proof.}   
Consider a locus $F_+$ of  
real cubic polynomials $z\mapsto z^3-3a^2z+b $ for which
the critical point $a$
is a preimage of the left fixed point (it is equivalent to $b=2a^3-2a$)
and $a<1/3$. 
By Branner and Douady  [BD] , there is a natural one-to-one correspondence
between $F_+$ and the 1/2-locus of quadratic polynomials
$z\mapsto z^2-c$ with $-2\leq c<-3/4.$ 
Hence, in $F_+$ there is only one Fibonacci
map (Theorem 1.1 ).
 On the other hand, conjugacy classes of cubic maps with escaping
critical point $a$ ( which means $b<2a^3-2a$)
are in one-to-one correspondence with $F+$ as well:  go toward the curve
$b=2a^3-2a$ along  external rays (this argument is due to Douady).\QED

From the last two lemmas we have immediate

{\bf Corollary 1. } All polynomial-like Fibonacci maps $f\in {\cal F}^-$ 
are quasi-symmetric\-ally conjugate.\medskip

{\bf Corollary 2.} Either all polynomial-like Fibonacci maps 
belong to the set ${\cal F}_0$ or  to its complement.\medskip

{\bf Proof.} For maps $f\in {\cal F}^-$ it follows from the last Corollary
and Lemma 5.3. For maps $f\in {\cal F}^+$  just observe that it belongs to
${\cal F}^0$ or its complement together with the renormalization.\QED

Now we will give an example of a polynomial-like map belongning to
${\cal F}_0$ which will yield that all Fibonacci polynomial-like maps
belong to ${\cal F}_0$. \medskip

{\bf Example.} Consider disjoint union of two intervals
$I=[-1,\lambda]$ and $J=[-c,-c+q\lambda^2]$ with positive $c, q, \lambda$,
$c$ is big, $\lambda$ is small. Let $f|I$ be a  quadratic
map $x\mapsto qx^2-c$, while $f|J$ be linear $x\mapsto \alpha x+b.$

Let us adjust parameters $\alpha,b,c,q,\lambda$ in such a way that 
$$0\mapsto -c\to -1\mapsto \lambda\mapsto -c+ q\lambda^2
\mapsto v\in [0,\lambda].$$
It yield the relations
$$q=c+\lambda\sim c ,\; \alpha={1+v\over (c+\lambda)\lambda^2}
\sim {1\over c\lambda^2},
     \;b=\alpha c-1\sim{1\over\lambda^2}\eqno (7.1)$$ 
It remains three free parameters $c, \lambda$ and $v$. Let us show that for
 $c^2\lambda^2<1$ this map is cubic-like. To this end consider  a disk 
$D=\{z: |z|<2\}$. On its boundary $\partial D$ our map acts as
$$f(z)=c(z^2-1)+\lambda z^2\sim c(z^2-1).$$
Hence,  $$3c<|f(z)|<5c \qquad {\rm for}\qquad z\in \partial D\eqno (7.2).$$ 
 Consider a disk $V=\{z: |z|<2c\}$
and its inverse image $U_1$ (under the quadratic map.) By (7.2), $U_1\subset D$
and $f: U_1\rightarrow V$ is a quadratic-like map. Moreover, $U_1\supset [-1,1]$ 
since $f[-1,1]=[-c,\lambda]\subset V.$

Furthermore, consider the preimage $U_2$ of $V$ under the linear map 
$z\to \alpha z+b$.
 It is a disk containing $J$  of radius
$$2c/\alpha\sim 2c^2\lambda^2<2$$
(by (7.1)). Hence, for big enough $c$ the closure of this disk is contained 
in $V$ and does not intersect cl$U_1$. So, $f: U_1\cup U_2\rightarrow V$
is a polynomial-like map.

Now one can adjust $v$ to get a Fibonacci map. Since $f$ has non-positive
Schwarzian derivative, it belongs to ${\cal F}^0$ 
provided $\lambda$ is sufficiently small (Lemma 5.2). \QED \medskip

 {\bf Renormalization of a quadratic-like Fibonacci map.}
  This procedure associate to a quadratic-like Fibonacci map
(of type (2))
 a cubic-like Fibonacci map (of type (2,1)).
It will complete the proof of Theorem 1.3 for quadratic-like Fibonacci maps
(in particular, for the quadratic polynomial).
We can restrict ourselves to the case of the quadratic
Fibonacci polynomial. Now let
us consider the beginning of the Yoccoz partition construction
(see [H] ).
 Draw a curve $S$ consisting of two external
rays through the fixed point $\alpha$ and an equipotential level $\gamma$.
We will obtain two pieces of level 0,
 $W^0$ (containing 0) and $W^0_1$ (containing $x_1$).
Define pieces  of level $n$  as
 $n$-fold preimages of the pieces of level 0. Denote by $W^n(x)$ the piece of
level $n$ containing $x$, set $W^n\equiv W^n(0)$.
Let us consider the piece $V\equiv W^4\supset T^4$ satisfying the property that 
$${\rm cl}\; W^4\subset W^3. \eqno (7.3)$$ 
Define a piece $U_1\equiv W^9\supset T^5$ as the pull-back of $V$
of order 5, and $U_2\supset J^5$ as the pull-back of $V$ of order 3.
One can check that  cl$U_1$ and cl$U_2$ are pairwise disjoint and are contained
in $V$ (it is a formal corollary from (7.3)). So, the map $g$ defined as
$f^5|U_1$ and $f^3|U_2$ is polynomial-like of type (2,1).\QED\medskip

{\bf Remark.} The above construction actually can be applied to any
non-infinitely renormalizable "persistently recurrent" quadratic polynomial
(see [L2]).

{\bf Geometry of $\omega(c)$ is not rigid.}
 We would like to show that parameter $a$
can really be changed in class ${\cal U}$, so the geometry of $\omega(c)$
is not rigid. The above Example provides us with a Fibonacci map of class
${\cal A}$ with arbitrary small $\lambda_0=1/c$. By Lemma 5.4, parameter
$a\asymp \lambda_0$ is getting arbitrary small as well. 
 Renormalizing $f$ if necessary we obtain a Fibonacci map
of class ${\cal A}^-$ with arbitrary small $a$. Now the surgery of \S 6
turns this map into a unimodal Fibonacci map with the same parameter $a$. 

{\bf Remark.} Actually, in order to vary parameter $a$ in class ${\cal A}$
 it is enough to observe that the renormalization turns $a$ into
$a/\root 3\of 2$.
\bigskip

\centerline{\bf \S 8. Polynomial-like property of analytic Fibonacci 
maps.}\medskip

In this section we will prove that analytic Fibonacci maps
$f\in{\cal E}$ become polynomial-like after apropriate 
renormalization.
Together with the results of the previous two sections it will complete the
proof of Theorem 1.3.

For an interval $I\subset {\bf R}$ denote by  
 by $D(I)$ the Eucledian disk based upon $I$ as the diameter.

{\QP{\bf Lemma 8.1} (see [S]). 
Let $\phi: P(I)\rightarrow P(J)$ be an analytic map which
 maps $I$ diffeomorphically onto $J$. Then $\phi D(I)\subset D(J).$\medskip}

{\bf Proof.} The interval $I$ is a Poincar\'{e} geodesic in $P(I)$, and
the disk $D(I)$ is its Poincar\'{e} neighborhood (of radius independent
of $I$). Since $\phi$  contracts the  Poincar\'{e} metric,
 we have the required. \QED

{\QP{\bf Lemma 8.2.} Let $f\in {\cal E}$ be an analytic Fibonacci map. Given $n$,
consider a disk $V=D(T_n)$ and its pull-backs $U_1\supset T_{n+1}$ and 
$U_2\supset J_{n+1}$ of 
order $u(n+2)$ and $u(n+1)$ correspondingly. Then cl$U_i$ are disjoint and are
contained in $V$.\medskip}

{\bf Proof.} Let $T_n=[[t_n,t_n']]$ with  $t_n$ being closer to $x_{u(n+2)}$.

The branch $\phi: V\rightarrow U_2$ of $f^{-u(n+1)}$ satisfies
the asumptions of Lemma 8.1, and hence $U_2\subset D(J_{n+1})$. By the same 
reason,  $fU_1\subset D(Q)$ where $Q\equiv [b,a]\ni x_1$ 
is the monotone pull-back of $T_n$ of order $u(n+2)-1$
( $b<x_1$ is the preimage of $t_n$ ).  

Now let  $X_{n-1}$ be the component of $T_{n-1}\ssm T_{n}$ adjacent to $t_n$.
 Since $\sum |X_n|<\infty$,  we can select such an $n$ that 
$$|X_n|<|X_{n-1}| \eqno (8-1).$$ 
By Lemma 4.1, the map $f^{u(n+2)-1}$ has a monotone continuation
beyond the point $b$ to the interval $W$ which is mapped onto $X_{n-1}$.
So, we have three interval map
 
$$f^{u(n+2)-1}: W\cup [b,x_1] \cup [x_1,a]\rightarrow
 X_{n-1}\cup [[t_n,x_{u(n+2)}]]\cup [[x_{u(n+2)},t_n']].\eqno (8-2)$$ 

Let $q=|x_{u(n+2)}|:|t_n|$. Applying the Schwarz lemma to (8-2)
taking into account (8-1) we get

$$\log {a-b\over a-x_1}\leq \log2 + \log{2\over 1+q},$$
so that
$${x_1-b\over a-x_1}\leq {3-q\over 1+q}.  \eqno (8-3)$$

Now let us take the $f$-preimage of $D(Q)$. Since $f^{-1}$ is just a square root
\break $\psi:\zeta\mapsto\sqrt{\zeta-x_1}$
on $D(Q)$, this preimage is contained in a domain based upon
$T_{n+1}$ with atitude  
 
$$h=|t_{n+1}| \sqrt{x_1-b\over a-x_1}\leq |t_{n+1}|\sqrt{3-q\over 1+q}\leq
   |t_{n+1}|/q<t_n.$$

Moreover, this domain is contained in the disk centered at zero of radius\break
$\max(t_{n+1},h)<t_n $. So, cl$U_1\subset V$.

 Let us show now that cl$U_1\cap{\rm cl}\,U_2=\emptyset$.\
If $a-x_1\geq x_1-b$ then $\psi D(Q)\subset D(T_{n+1})$,
and the statement follows.
Assume that $x_1-b > a-x_1$. Then one can check the following elementary fact
about the square root map:
$\psi D[b,a]$ is convex if and only if $x_1-b\leq 3(a-x_1)$. By (8-3), the
last estimate holds, so $\psi D(Q)$ is convex. Hence, 
$\psi D(Q)\cap D(J_{n+1})=\emptyset$, and we are done.
\QED
\bigskip

\centerline{\bf Appendix. Schwarz Lemma and Koebe Principle.} \medskip

We refer the reader to [Y], [G2], [Sw1-2], [MS] and [S] 
for the following technical background.
 
Let us consider four points $a<b<c<d$ and two nested intervals $L=[a,d]$
and $H=[b,c]$. The {\sl Poincar\'{e} length} of $H$ in $L$ is the logarithm
of an appropriate cross-ratio:
$$[H:L]=\log{(d-b)(c-a)\over (d-c)(b-a)}.$$
Let $g: (L,H)\rightarrow (L',H')$ be a $C^3$ diffeomorphism,
$$Sg={g'''\over g'}-{3\over 2}\left({g''\over g'}\right)^2$$
 be its Schwarzian derivative.

{\QP{\bf Schwarz Lemma.} If $g$ has non-negative Schwarzian derivative then 
it contracts Poincar\'{e} length:
$$[H':L']\leq [H:L].$$\medskip}

{\QP{\bf Koebe Principle.} Let $g$ has non-negative Schwarzian derivative.
If $[H:L]\leq\ell$ then
$$\left |{g'(x)\over g'(y)}\right |\leq K(\ell)$$
for any $x,y\in H$. Moreover, $K(\ell)=1+O(\ell)$ as $\ell\to 0$.\medskip} 

One can essentially extend the range of applications of these results 
combining the  Schwarzian derivative condition on some intervals with
bounded non-linearity on others. Let us consider a chain of (closed) interval
diffeomorphisms
$$I_1\rightarrow J_1\rightarrow...\rightarrow I_n\rightarrow J_n$$
where $g_i : I_i\rightarrow J_i$ have negative Schwarzian derivative
while $h_i : J_i\rightarrow I_{i+1}$ are just $C^2$ smooth.
Set $F=h_n\circ g_n\circ...\circ h_1\circ g_1$. Let $G_i\subset {\rm int} I_i$ and
$H_i\subset {\rm int} J_i$ be closed subintervals related by diffeomorphisms.

Denote by ${\bf h}$ the family of maps $h_i$, by ${\bf I}$ the family of intervals
$I_i$ etc.
Let\break $\|h_i\|=\max|h''(x)/h'(x)|$,
 $\|{\bf h}\|=\max \|h_i\| $ 
be the {\sl maximal nonlinearity } of ${\bf h}$, $|{\bf I}|=\sum |I_i|$ be the
{\sl total length} of ${\bf I}$.

{\QP{\bf Schwarz Lemma (smooth version).} Expantion of the Poincar\'{e} length
by the map $F$ is controlled by ${\bf h}$ in the following manner
$$[H_n:J_n]\leq [G_n:I_n]+O(|{\bf J}|)$$
with the constant depending on ${\|\bf h\|}$.\medskip}

{\QP{\bf Koebe Principle (smooth version).} Distortion of $F|G_1$ can 
be estimated as follows:

$$\left |{F'(x)\over F'(y)}\right |\leq K(\ell; |{\bf h}|, |{\bf J}|)$$
where $K=1+O(\ell+|{\bf J}|)$ as  $|{\bf J}|, \ell\to 0$
with the constant depending on $|{\bf h}|$.

\bigskip
\centerline{\bf References.}\smallskip

{\item [BD]} B.Branner \& A.Douady. Surgery on complex polynomials, Preprint
    Matematic Institut, Denmark, 1987-05. 
{\item [BH]} B.Branner \& J.H.Hubbard. The iteration of cubic polynomials, Part II   : patterns and parapatterns, Acta Math., to appear.\smallskip
{\item [BL1]} A.Blokh \& M.Lyubich. Non-existence of wandering intervals and
    structure of topological attractors for one dimensional dynamical systems,
    Erg. Th. \& Dyn Syst. {\bf 9} (1989), 751-758.\smallskip
 {\item [BL2]} A.Blokh \& M.Lyubich. Measurable dynamics of S-unimodal maps
    of the interval, Preprint IMS Stony Brook. \# 1990/2\smallskip
{\item [BL3]} A.Blokh \& M.Lyubich. Measure and dimension of solenoidal
   attractors of one-dimensional dynamical systems, Comm. Math. Phys.,
    {\bf 127} (1990), 573-583.\smallskip
{\item [DH]} A.Douady \& J.H.Hubbard. On the dynamics of polynomial-like maps,
    Ann. Sc. Ec. Norm. Sup. {\bf 18} (1985), 287-343.\smallskip
{\item [G1]} J. Guckenheimer. Sensitive dependence to initial conditions for
   one-dimensional maps, Comm. Math. Phys. {\bf 70} (1979), 133-160.\smallskip
{\item [G2]} J.Guckenheimer. Limit sets of S-unimodal maps with zero entropy,
   Comm. Math. Phys., {\bf 110}, 655-659.\smallskip
{\item [H]} J.H.Hubbard, according to J.-C.Yoccoz. Puzzles and quadratic
tableaux. Preprint, 1990.
{\item [HK]} F. Hofbauer and G. Keller. Some remarks on recent results 
about S-unimodal maps. Preprint, 1990.
{\item [L1]} M.Lyubich. Non-existence of wandering intervals and
    structure of topological attractors for one dimensional dynamical systems,
    Erg. Th. \& Dyn Syst. {\bf 9} (1989), 737-750.\smallskip
{\item [L2]} M.Lyubich.On the Lebesgue measure of the Julia set of a 
     quadratic polynomial, Preprint IMS, 1991/10. \smallskip
{\item [M]} J.Milnor. On the concept of attractor. Comm. Math. Phys,
     {\bf 99} (1985), 177-195, and {\bf 102} (1985), 517-519. \smallskip
{\item [MT]} J.Milnor \& W.Thurston. On iterated maps of the interval,
pp. 465-563 of ``Dynamical Systems, Proc. U. Md., 1986-87,
ed. J. Alexander, Lect. Notes Math., 1342, Springer 1988.\smallskip
{\item [Ma]} M.Martens. Cantor attractors of unimodal maps. Preprint, 1990.
\smallskip
{\item [MS]} W.de Melo \& S. van Strien. A structure theorem in one-dimensional
dynamics. Preprint, 1986.\smallskip
{\item [MMSS]} M.Martens \& W.de Melo \& S. van Strien \& D.Sullivan.
Bounded geometry and measure of the attracting Cantor set of
quadratic-like maps. Preprint, 1988.\smallskip
{\item [NS]} T.Nowicki \& S. van Strien. Invariant measures exist under a
summability condition for unimodal maps. Preprint.\smallskip
{\item [PTT]} I.Procaccia \& S.Thomae \& C.Tresser. First return maps as
a unified renormalization scheme for dynamical systems.
Physical Review A, {\bf 35} (1987
), n 4, 1884-1900.\smallskip
{\item [S]} D.Sullivan. Bounds, quadratic differentials, and renormalization
conjectures. Preprint,  1990. \smallskip
{\item [Sh]} K.Shibayama. Fibonacci sequence of stable periodic orbits
for one-parameter families of $C^1$-unimodal mappings. Preprint.\smallskip
{\item [Sw1]} G. Swiatek. Rational rotation numbers  for maps of the circle.
Comm. Math. Phys., {\bf 119} (1988), 109-128.\smallskip 
{\item [Sw2]} G. Swiatek. Bounded distortion properties of one-dimensional maps.
Preprint SUNY, Stony Brook, 1990/10. \smallskip
{\item [TV]} F.Tangerman \& P.Veerman.  Scaling in circle maps I.
 Preprint, SUNY, Stony Brook, 1990/8. \smallskip
[V] A. Vershik. Generalized notations and dynamical systems. Preprint, 1991.
   \smallskip
[Y1] J.-C. Yoccoz. Il n'y a pas de contre-exemple de Denjoy analytiques.
C.R.Acad. Sci. Paris, {\bf 289} (1984), 141-144. \smallskip
[Y2] J.-C. Yoccoz. Manuscript, 1990. \smallskip

\end